\def\draft{n}
\documentclass{amsart}
\usepackage{fullpage,amssymb,epic,eepic,epsfig,amscd,pb-diagram,lamsarrow,pb-lams}

\theoremstyle{plain}

\newtheorem{theorem}{Theorem}
\newtheorem{proposition}{Proposition}[section]
\newtheorem{lemma}[proposition]{Lemma}
\newtheorem{corollary}[proposition]{Corollary}

\theoremstyle{definition}
\newtheorem{definition}[proposition]{Definition}

\theoremstyle{remark}

\newtheorem{remark}[proposition]{Remark}

\def\printname#1{
	\if\draft y
		\smash{\makebox[0pt]{\hspace{-0.5in}
			\raisebox{8pt}{\tt\tiny #1}}}
	\fi
}

\newcommand{\psdraw}[2]
         {\begin{array}{c} \hspace{-1.3mm}
	\raisebox{-4pt}{\epsfig{figure=draws/#1.eps,width=#2}}
	\hspace{-1.9mm}\end{array}}

\newlength{\standardunitlength}
\catcode`\@=11
\long\def\@makecaption#1#2{%
     \vskip 10pt

\setbox\@tempboxa\hbox{%\ifvoid\tinybox\else\box\tinybox\fi
       \small\sf{\bfcaptionfont #1. }\ignorespaces #2}%
     \ifdim \wd\@tempboxa >\captionwidth {%
         \rightskip=\@captionmargin\leftskip=\@captionmargin
         \unhbox\@tempboxa\par}%
       \else
         \hbox to\hsize{\hfil\box\@tempboxa\hfil}%
     \fi}
\font\bfcaptionfont=cmssbx10 scaled \magstephalf
\newdimen\@captionmargin\@captionmargin=2\parindent
\newdimen\captionwidth\captionwidth=\hsize
\catcode`\@=12

\newcommand{\tr}{\operatorname{tr}}

\def\lbl#1{\label{#1}\printname{#1}}

\newcommand{\eatline}{\vspace{-\baselineskip}}

\def\BZ{\mathbb Z}
\def\BQ{\mathbb Q}
\def\BR{\mathbb R}
\def\BC{\mathbb C}

\def\A{\mathcal A}
\def\B{\mathcal B}

\def\cD{\mathcal D}
\def\D{\Delta}

\def\a{\alpha}

\def\La{\Lambda}
\def\Lz{\Lambda_{\BZ}}
\def\lb{\lambda}
\def\Ga{\Gamma}
\def\S{\Sigma}
\def\Sp{\Sigma^p}
\def\s{\sigma}
\def\vars{\varsigma}

\def\ihs{integral homology 3-sphere}
\def\qhs{rational homology 3-sphere}

\def\fti{finite type invariant}

\def\la{\langle}
\def\ra{\rangle}

\def\w{\omega}

\def\g{\gamma}

\def\e{\epsilon}
\def\Ga{\Gamma}
\def\d{\delta}
\def\b{\beta}

\def\Th{\Theta}
\def\s{\sigma}

\def\lk{{\text{lk}}}

\newcommand{\Ker}{\operatorname{Ker}}

\def\Sei{\mathrm{Sei}}

\def\pt{\mathbb Z[\![t_1 ,\dots,t_m ]\!]}
\def\sminus{\smallsetminus}

\def\ti{\widetilde}

\def\wheel{\operatorname{\psfig{figure=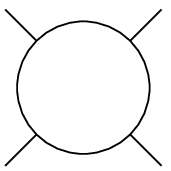,height=0.1in}}}
\def\ostar{\circledast}

\def\tline{\uparrow}

\def\strutb#1#2#3{\overset{#1}{\underset{#2}{ 
\begin{array}{c} \vspace{0.0cm}
\uparrow 
\vspace{-0.25cm} \\        %\vspace{-0.25cm}
| \vspace{-0.45cm} \\      %\vspace{-0.45cm}
\bullet \vspace{0.00cm}   %\vspace{0.00cm}
\end{array} }}\! #3}

\def\st#1#2#3{\overset{#2}{\underset{#1}\uparrow} #3}

\def\AS{\mathrm{AS}}
\def\IHX{\mathrm{IHX}}
\def\GA{\A^{\mathrm{gp}}}
\def\GAz{\A^{\mathrm{gp},0}}

\def\clover{clover}

\def\Lloc{\Lambda_{\mathrm{loc}}}

\def\loc{\mathrm{loc}}

\def\Zrat{Z^{\mathrm{rat}}}
\def\Zratc{\check{Z}^{\mathrm{rat}}}
\def\Zratcw{\check{Z}^{\mathrm{rat},\wheel}}

\def\Zratw{Z^{\mathrm{rat},\wheel}}
\def\Zc{\check{Z}}
\def\Zcw{\check{Z}^{\wheel}}
\def\Zh{\hat{Z}}
\def\Z{Z}
\def\Zw{Z^{\wheel}}

\def\intrat{\int^{\mathrm{rat}}}

\def\twist{\tau}                  
\def\twistrat{\tau^{\mathrm{rat}}}  
\def\Res{\mathrm{Res}}
\def\lift{\mathrm{Lift}}
\def\liftrat{\mathrm{Lift}^{\mathrm{rat}}}

\def\sig{$\sigma$ignature}
\def\ch{\mathrm{ch}}

\def\Aut{\mathrm{Aut}}
\def\hair{\mathrm{Hair}}
\def\hairnu{\mathrm{Hair}^{\Om}}
\def\phia{\phi_{k\to k+h}}
\def\phim{\phi_{t\to t e^h}}
\def\phimp{\phi_{t\to t e^{h'}}}
\def\phimm{\phi_{t\to t e^h e^{h'}}}

\def\Herm{\mathrm{Herm}}
\def\Sym{\mathrm{Sym}}
\def\eyes{\operatorname{\psfig{figure=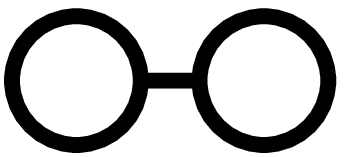,height=0.1in}}}
\def\Om{\Omega}
\def\homega{\widehat\Omega}
\def\bgp{\beta^{\mathrm{gp}}}
\def\kbr{K_{\mathrm{br}}}
\def\conh{\mathrm{con}_{\{h\}}}

\def\pt{\partial}
\def\pr{\mathrm{pr}}
\def\longto{\longrightarrow}

\begin{document}

%%%%%%%%%%%%%%%%%%%%%%{page1}

\title[Finite type invariants of cyclic branched covers]{
Finite type invariants of cyclic branched covers}

\author{Stavros Garoufalidis}
\address{School of Mathematics \\
          Georgia Institute of Technology \\
          Atlanta, GA 30332-0160, USA. }
\email{stavros@math.gatech.edu}
\author{Andrew Kricker}
\address{Institute of Mathematics, Hebrew University at Giv'at-Ram, Jerusalem
   91904, Israel.}
\email{kricker@math.huji.ac.il}

\thanks{S.G. was partially supported %by an NSF grant DMS-98-00703 and 
by an Israel-US BSF grant, and A.K. was partially
supported by a JSPS Fellowship.
        This and related preprints can also be obtained at
{\tt http://www.math.gatech.edu/$\sim$stavros } 
%and {\tt http: }
\newline
1991 {\em Mathematics Classification.} Primary 57N10. Secondary 57M25.
}

\date{
This edition: July 1, 2002 \hspace{0.5cm} First edition: August 8, 2000.}

%\dedicatory{Preliminary notes. Please do not distribute under
%any circumstances!}

\begin{abstract}
Given a knot in an integer homology sphere, one can construct a family of
closed 3-manifolds (parametrized by the positive integers), namely the cyclic 
branched coverings of the knot. In this paper we give a formula for the 
the Casson-Walker invariants of these 3-manifolds in terms of 
residues of a rational function (which measures the 2-loop part of the
Kontsevich integral of a knot) and the signature function of the knot.
Our main result actually computes the LMO invariant of cyclic branched covers
in terms of a rational invariant of the knot and its signature function.
\end{abstract}

\maketitle

%\tableofcontents

%%%%%%%%%%%%%%%%% the text file

\section{Introduction}
\lbl{sec.intro}
\subsection{History} 
\lbl{sub.history}

One of the best known integer-valued concordance invariants of a knot $K$
in an integer homology sphere $M$ is its (suitably normalized)
{\em \sig\ function} $\s(M,K):S^1\to\BZ$ defined for all complex numbers
of absolute value $1$, see for example \cite{Ka}. The \sig\ function
and its values at complex roots of unity are closely 
related to a sequence (indexed by a natural number $p$, not necessarily prime)
of closed 3-manifolds, the {\em $p$-fold 
cyclic branched coverings} $\Sp_{(M,K)}$,
associated to the pair $(M,K)$ and play
a key role in the approach to knot theory via surgery theory.

It is an old problem to find a formula for the Casson-Walker invariant
of cyclic branched covers of a knot. For two-fold branched covers,
 Mullins used skein theory of the Jones polynomial to 
show that for all knots $K$ in $S^3$ such that $\S^2_{(S^3,K)}$ is a \qhs ,
there is a linear relation between $\lb(\S^2_{(S^3,K)})$, $\s_{-1}(S^3,K)$
and the logarithmic derivative of the Jones polynomial of $K$ at $-1$, 
\cite{Mu}. A different approach was taken by the first author in \cite{Ga},
where the above mentioned linear relation was deduced and explained from the 
wider context of \fti s of knots and 3-manifolds.    

For $p>2$, Hoste, Davidow and Ishibe 
studied a partial case of the above problem for Whitehead doubles of 
knots, \cite{Da,Ho,I}. 

However, a general formula was missing for $p >2$. Since the map
$(M,K) \to \lb(\Sp_{(M,K)})$ is not a concordance invariant of 
$(M,K)$, it follows that a formula for the Casson invariant of cyclic
branched coverings should involve more than just the 
{\em total $p$-\sig } $\s_p$ (that is, $\sum_{\w^p=1}\s_\w$).

In \cite{GR}, a conjecture for the Casson invariant of cyclic
branched coverings was formulated. The conjecture involved the 
total signature and the sums over complex roots of unity, of a {\em rational
function} associated to a knot. The rational function in question 
was the 2-loop part of a rational lift $\Zrat$ of the Kontsevich integral of 
a knot.

In \cite{GK2} the authors constructed this rational lift, combining the 
so-called {\em surgery view of knots} (see \cite{GK1}) with the full aparatus
of perturbative field theory, formulated by the Aarhus integral and its
function-theory properties.

The goal of the present paper is to prove the missing formula of the
Casson invariant of cyclic branched coverings, under the mild assumption
that these are rational homology spheres. In fact, our methods will give
a formula for the LMO invariant of cyclic branched coverings in terms of
the \sig\ function and residues of the $\Zrat$ invariant.
 
Our main Theorem \ref{thm.1} will follow from a formal calculation, 
presented in Section \ref{sub.formal}. This illustrates the relation between
the formal properties of the $\Zrat$ invariant and the geometry of the
cyclic branched coverings of a knot.

\subsection{Statement of the results}
\lbl{sub.results}

Let us call a knot $(M,K)$ {\em $p$-regular} iff $\Sp_{(M,K)}$ is a \qhs .
We will call a knot {\em regular} iff it is $p$-regular for all $p$.
It is well-known that $(M,K)$ is $p$-regular iff its Alexander polynomial 
$\D(M,K)$ has no complex $p$th roots of unity.

Let $Z$ denote the LMO invariant of a knot (reviewed in Section 
Section \ref{sec.review}),  and let $\twistrat$ denote the the twisting map 
of Definition \ref{def.twistrat} and $\lift_p$ denote the lifting map of
Section \ref{sub.deflift}.

\begin{theorem}
\lbl{thm.1}
For all $p$ and $p$-regular pairs $(M,K)$ we have
$$\Z(\Sp_{(M,K)})=e^{\s_p(M,K)\Theta/16} 
\lift_p \circ \twistrat_{\a_p}  \circ \Zrat(M,K)  \in  
\A(\phi).
$$
where $\a_p=\nu^{-(p-1)/p} \in \A(\ostar)$, $\nu=Z(S^3,\mathrm{unknot})$.
\end{theorem}

The proof of Theorem \ref{thm.1} is a formal computation, given in Section
\ref{sec.reduction}, that involves the 
rational invariant $\Zrat$ and its function-theory properties, phrased in 
terms of operations (such as twisting and lifting) on diagrams. 
In a sense, the $\Zrat$ invariant is defined by using properties of the
universal abelian cover of knot complements. Since the universal abelian
cover maps onto every cyclic branched cover, it is not too surprising
that the $\Zrat$ invariant appears in a formula for the LMO invariant of cyclic
branched covers. The presence of the signature function is
a framing defect of the branched covers. It arises because we need to
normalize the 3-manifold invariants by their values at a $\pm 1$-framed unknot.
These values are some universal constants, whose ratio (for positively versus
negatively unit framed unknot) is given by the signature term.
We do not know of a physics explanation of the above formula in terms
of anomalies.

\begin{remark}
\lbl{rem.ask}
Twisting and lifting are important operations on diagrams with beads 
that commute with the operation of integration, see Propositions 
\ref{prop.ztwist} and \ref{prop.liftint}. For properties of 
the twisting operation, see Lemma \ref{lem.twistrat} in Section 
\ref{sec.twist}. For a relation between our notion of twisting and the notion 
of {\em wheeling} (introduced in \cite{A0} 
and studied in \cite{BLT,BL}), see Section \ref{sec.wheel}. For properties
of the lifting operation, see Section \ref{sec.lifting}.
Twisting and Lifting 
are closely related to the Magic Formula for the Kontsevich integral of
a Long Hopf Link \cite{BLT}, and to rational framings, \cite{BL}.
\end{remark}

The next corollary gives a precise answer for the value of the Casson-Walker
invariant of cyclic branched covers as well as its growth rate (as $p \lim
\infty $), in terms of the $2$-loop part of the Kontsevich integral and the
\sig\ function. In a sense, the \sig\ function and the
2-loop part of the Kontsevich integral are generating function for the
values of the Casson-Walker invariant of cyclic branched covers.

\begin{corollary}
\lbl{cor.1}
$\mathrm{(a)}$ For all $p$ and $(M,K)$ and $p$-regular, we have
$$
\lb (\S^p_{(M,K)})= \frac{1}{3} \, \Res^{t_1,t_2,t_3}_p Q(M,K)(t_1,t_2,t_3)
+ \frac{1}{8} \, \s_p(M,K).
$$
Note the difference between the normalization of $\Res_p$ of 
\cite[Section 1.5]{GR} and that of Section \ref{sub.connection}. \newline
$\mathrm{(b)}$
For all regular pairs $(M,K)$, we have
$$
\lim_{p \to \infty} \frac{\lb (\S^p_{(M,K)})}{p} = \frac{1}{3}
\int_{S^1 \times S^1} Q(M,K)(s) d\mu(s) +
\frac{1}{8} \int \s_s(M,K) d\mu(s)
$$
where $d\mu$ is the Haar measure.
\end{corollary}

In other words, the Casson invariant of cyclic branched coverings grows
linearly with respect to the degree of the covering, and the growth rate
is given by the average of the $Q$ function on a torus and the {\em
total \sig } of the knot (i.e., the term $\int \s_s(M,K) d\mu(s)$ above).
The reader may compare this with the following theorem
of Fox-Milnor, \cite{FM} 
which computes the torsion of the first homology of cyclic branched covers 
in terms of the Alexander polynomial, and the growth rate of it in terms
of the {\em Mahler measure} of the Alexander polynomial: 

\begin{theorem}\cite{FM}
$\mathrm{(a)}$
Let $\b_p(M,K)$ denote the order of the torsion
subgroup of $H_1(\S^p_{(M,K)},\BZ)$.  Assuming that $(M,K)$ is $p$-regular,
we have that:
$$
\b_p(M,K)=\prod_{\w^p=1}| \D(M,K)(\w)|.
$$
$\mathrm{(b)}$
If $(M,K)$ is regular, it follows that
$$
\lim_{p \to \infty} \frac{\log \b_p(M,K)}{p}
= \int_{S^1} \log (|\D(M,K)(s)|) d\mu(s).
$$
\end{theorem}

In case $(M,K)$ is not regular, $\mathrm{(b)}$ still holds, as was shown
by Silver and Williams, \cite{SW}.  

\subsection{Plan of the proof}
\lbl{sub.plan}
In Section \ref{sec.reduction}, we review the definition of $\Zrat$, and
we reduce Theorem \ref{thm.1} to Theorem \ref{thm.2} (which concerns
signatures of surgery presentations of knots) and Theorem \ref{thm.3}
(which concerns the behavior of the $\Zrat$ invariant under coverings 
of knots in solic tori).

Section \ref{sec.comparison} consists entirely of topological facts about
the surgery view of knots, and shows Theorem \ref{thm.2}.

Sections \ref{sec.twist} and \ref{sec.lifting} introduce the notion of
twisting and lifting of diagrams, and study how they interact with the
formal diagrammatic properties of the $\Zrat$ invariant. As a result,
we give a proof of Theorem \ref{thm.3}.

In Section \ref{sub.deg2} we prove Corollary \ref{cor.1}.

Finally, we give two alternative versions of Theorem \ref{thm.1}: 
in Section \ref{sec.remember} in terms of an invariant of branched covers
that remembers a lift of the knot, and in Section \ref{sec.wheel} in terms of 
the wheeled rational invariant $\Zratw$.
 
\subsection{Recommended reading}
The present paper uses at several points a simplified version of the notation
and the results of \cite{GK2} presented for knots rather than boundary links.
Therefore, it is a good idea to have a copy of \cite{GK2} available.

\subsection{Acknowledgement}
We wish to thank L. Rozansky and D. Thurston and especially 
J. Levine and T. Ohtsuki for stimulating conversations and their 
support. The first author was supported
by an Israel-US BSF grant and the second author was supported by a JSPS
Fellowship.

\tableofcontents

\section{A  reduction of Theorem \ref{thm.1}}
\lbl{sec.reduction}

In this section we will reduce Theorem \ref{thm.1} to two theorems;
one involving properties of the invariant $\Zrat$ under lifting and
integrating, and another involving properties of the 
\sig\ function. Each will be dealt with in a subsequent section.

\subsection{A brief review of the rational invariant $\Zrat$}
\lbl{sec.review}

In this section we briefly explain where the rational invariant takes
values and how it is defined. 
The invariant $\Zrat(M,K)$ is closely related to the surgery view of
pairs $(M,K)$ and is defined in several steps explained in \cite{GK1} and 
below, with some simplifications since we will be dealing {\em exclusively}
with knots and not with 
boundary links, \cite[Remark 1.6]{GK2}. In that case, 
the rational invariant $\Zrat$ takes values in the subset
$$
\GAz(\Lloc)=\B \times \GA(\Lloc) \hspace{1.5cm} \text{ of }
\hspace{1.5cm} \A^0(\Lloc)=\B \times \A(\Lloc)
$$
where 
\newline
$\bullet$ $\Lz=\BZ[\BZ]=\BZ[t^{\pm 1}]$, 
$\La=\BQ[\BZ]=\BQ[t^{\pm 1}]$ and $\Lloc=
\{p(t)/q(t), p, $ $q \in \BQ[t^{\pm 1}], q(1)=\pm 1\}$, the  
localization of $\La$ with respect to the multiplicative set
of all Laurent polynomials of $t$ that evaluate to $1$ at $t=1$. For future
reference, $\La$ and $\Lloc$ are rings with involution $t\leftrightarrow 
t^{-1}$,
selected group of units $\{ t^n \, | \, n \in \BZ\}$ and ring homomorphisms
to $\BZ$ given by evaluation at $t=1$.
\newline
$\bullet$
$\Herm(\Lz\to\BZ)$ is the set of {\em Hermitian matrices} $A$
over $\Lz$, invertible over $\BZ$, and $B(\Lz\to\BZ)$ ({\em abbreviated by
$\B$}) denote the quotient 
of $\Herm(\Lz\to\BZ)$ 
modulo the equivalence relation generated by the move: 
$A \sim B$ iff $A\oplus E_1=
P^\ast(B\oplus E_2)P$, where $E_i$ are diagonal matrices with $\pm 1$
on the diagonal and $P$ is either an elementary matrix (i.e., one that differs
from the diagonal at a single nondiagonal entry)
or a diagonal matrix with monomials in $t$ in the diagonal.
\newline
$\bullet$
$\A(\Lloc)$ is the (completed) graded algebra over $\BQ$ spanned by
 trivalent graphs
(with vertex and edge orientations) whose edges are labeled by elements
in $\Lloc$, modulo the $\AS,\IHX$ relations and the Multilinear and Holonomy
relations of \cite[Figures 1,2, Section 3]{GK2}.
The degree of a graph is the number of its trivalent vertices and the
multiplication of graphs is given by their disjoint union.
$\GA(\Lloc)$ is the set of {\em group-like} elements of $\A(\Lloc)$, that is
elements of the form $\exp(c)$ for a series $c$ of connected graphs.

So far, we have explained where $\Zrat$ takes values. In order to recall the 
definition of $\Zrat$, we need to consider unitrivalent graphs as well
and a resulting set $\GA(\ostar_X,\La)$ explained in detail in Section
\ref{sec.twist}. Then, we proceed as follows:
\newline
$\bullet$
Choose a {\em surgery presentation} $L$ for $(M,K)$, that is a null homotopic 
framed link $L$ (in the sense that each component of $L$ is a 
null homotopic curve in $ST$) in a standard solid torus $ST \subset S^3$ such 
that its linking matrix is invertible over $\BZ$ and such that
$ST_L$ can be identified with the complement of a tubular neighborhood 
of $K$ in $M$. 
\newline
$\bullet$
Define an invariant $\Zratc(L)$ with values in 
$\GA(\ostar_X,\La)$ where $X$ is a set in 1-1 correspondence with the 
components of $L$.
\newline
$\bullet$
Define an integration $\intrat \, dX: \GA(\ostar_X,\La)\to
\GAz(\Lloc)$ as follows. Consider an integrable element $s$, that is one 
of the form
\begin{equation}
\lbl{eq.Gaussd}
s = \exp_\sqcup \left( \frac{1}{2} \sum_{i,j}
 \strutb{x_j}{x_i}{M_{ij}} \right) \sqcup R,
\end{equation}
with $R$ a series of $X$-substantial diagrams (i.e., diagrams that do not
contain a strut component). Notice that $M$, the
{\em covariance matrix} of $s$, and $R$, the $X$-{\em substantial 
part} of $s$, are uniquely determined by $s$, and define
$$
\intrat dX(s) = \left(M,
\left\la  
\exp_{\sqcup}\left( -\frac{1}{2} \sum_{i,j} 
\strutb{x_j}{x_i}{M_{ij}^{-1}} \right)
,R\right\ra_{X} \right).
$$
In words, $\intrat$-integration is gluing the legs of the $X$-substantial 
graphs in $X$ using the negative inverse covariance matrix. 
\newline
$\bullet$
Finally, define 
\begin{equation}
\lbl{eq.zdef}
\Zrat(M,K) =
\frac
{  \intrat dX \,\, \Zratc(L)  }
{c_+^{\sigma_+(B)} c_-^{\sigma_-(B)}} \in \GAz(\Lloc),
\end{equation}
where $c_\pm=\int dU \Zc (S^3,U_{\pm})$ are some universal constants
of the unit-framed unknot $U_{\pm}$.

The following list of frequently asked questions may motivate a bit the
construction of $\Zrat$:
\newline
{\em Question:}
Why is $\Zrat(M,K)$ an invariant of $(M,K)$ rather than of $L$? 
\newline
{\em Answer:} Because
for fixed $(M,K)$ any two choices for $L$ are related by a sequence of
Kirby moves, as shown by the authors in \cite[Theorem 1.1]{GK2}. 
Even though $\Zratc(L) \in \GA(\star_X,\La)$ is not invariant under Kirby 
moves, it becomes so after $\intrat$-integration. 
\newline
{\em Question:} 
Why do we need to introduce the $\intrat$-integration? 
\newline
{\em Answer:} To make $L\to\Zratc(L)$ invariant under Kirby moves on $L$.
\newline
{\em Question:}
Why do we need to consider diagrams with beads in $\Lloc$?
\newline
{\em Answer:} Because $\intrat$-integration glues struts
by inverting the covariance matrix $W$. If $W$ is a Hermitian
matrix over $\La$ which is invertible over $\BZ$, after $\intrat$-integration
appear diagrams with beads entries of $W^{-1}$, a matrix defined over $\Lloc$.

\begin{remark}
\lbl{rem.normalizations}
$Z$ stands for the Kontsevich integral of framed links
in $S^3$, extended to an invariant of links in 3-manifolds by 
Le-Murakami-Ohtsuki, \cite{LMO}, and identified with the Aarhus integral
in the case of links in \qhs s, \cite[Part III]{A}. In this paper we will
use exclusively the Aarhus integral $\int$ and its rational generalization
$\intrat$, whose properties are closely related to
function-theoretic properties of functions on Lie groups and Lie algebras.

By convention, $\Zrat$ contains {\em no} wheels and no $\Om$ terms.
That is, $\Zrat(S^3,U)=1$. On the other hand, $\Z(S^3,U)=\Om$.
Note that $\Zratc(L)$ equals to the connect sum of copies of $\Om$ (one
to each component of $L$) to $\Zrat(L)$. 
\end{remark}

\subsection{Surgery presentations of cyclic branched covers}
\lbl{sub.surgery}

Fix a surgery presentation $L$ of a pair $(M,K)$. 
We begin by giving a surgery presentation of $\Sp_{(M,K)}$. 
Let $L^{(p)}$ denote the preimage of $L$ under the $p$-fold cover $ST\to ST$.
It is well-known that $L^{(p)}$ can be given a suitable framing so that
$\Sp_{(M,K)}$ can be identified with $S^3_{L^{(p)}}$, see \cite{CG}.

It turns out that the total $p$-\sig\ can be calculated from the
the linking matrix of the link $L^{(p)}$. In order to state the result,
we need some preliminary definitions.
For a symmetric matrix $A$ over $\BR$, let $\s_+(A), \s_-(A)$ 
denote the number of positive
and negative eigenvalues of $A$, and let $\s(A),\mu(A)$ denote the {\em
signature} and {\em size} of $A$. Obviously, for nonsingular $A$, we have
$\s(A)=\s_+(A)-\s_-(A)$ and $\mu(A)=\s_+(A)+\s_-(A)$.

Let $B$ (resp. $B^{(p)}$) denote the linking matrix of the framed link 
$L$ (resp. $L^{(p)}$) in $S^3$. We will show later that 
\begin{theorem}
\lbl{thm.2}(Proof in Section \ref{sub.Yview})
With the above notation, we have
$$
\s_p(M,K)=\s(B^{(p)})-p\s(B) \,\,\, \text{ and } \,\,\, \mu(B^{(p)}) = p\mu(B).
$$
\end{theorem}

\subsection{A formal calculation}
\lbl{sub.formal}

Assuming the existence of a suitable maps $\lift_p$ and $\twistrat$, take
residues of Equation \eqref{eq.zdef}. We obtain that

\begin{xalignat*}{2}
\lift_p \circ \twistrat_\a \circ \Zrat (M,K)  &= 
\lift_p \circ \twistrat_\a \left(
\frac{ 
\intrat dX \,\, \Zratc(L)}{c_+^{\sigma_+(B)} c_-^{\sigma_-(B)}}
\right) & \\
&= \lift_p \left(
\frac{ 
\intrat dX \,\, \twistrat_\a \Zratc(L)}{c_+^{\sigma_+(B)} c_-^{\sigma_-(B)}}
\right) & \text{by Theorem \ref{prop.ztwist}} \\
&=  
\frac{
\lift_p \left(\intrat dX \,\, \twistrat_\a \Zratc(L)  \right)}
{c_+^{p\sigma_+(B)} c_-^{p\sigma_-(B)}} & \text{by Remark \ref{rem.noskeleton}}
\\
&=
\frac{
\lift_p \left(\intrat dX \,\, \twistrat_\a \Zratc(L)  \right)}
{\left(\sqrt{\frac{c_+}{c_-}}\right)^{p\sigma(B)} 
{\left(\sqrt{c_+c_-}\right)^{p\mu(B)}}}. &
\end{xalignat*}

\noindent
Adding to the above the term corresponding to the total 
$p$-\sig\ $\s_p(M,K)$ of $(M,K)$, and using 
the identity $c_+/c_-=e^{-\Theta/8}$ (see \cite[Equation (19), 
Section 3.4]{BL}) it follows that

\begin{xalignat*}{2}
\lift_p \circ \twistrat_{\a_p} \circ \Zrat(M,K) e^{\s_p(M,K)\Theta/16} 
&= 
\lift_p \circ \twistrat_{\a_p} \circ \Zrat(M,K) 
\left(\sqrt{\frac{c_+}{c_-}}\right)^{-\sigma_p(M,K)} 
& %\text{by \cite[Equation 19, Section 3.4]{BL}} 
\\ &= 
\frac{
\lift_p \left(\intrat dX \,\, \twistrat_{\a_p} \Zratc(L)  \right)}
{\left(\sqrt{\frac{c_+}{c_-}}\right)^{\sigma( B^{(p)} ) } 
{\left(\sqrt{c_+c_-}\right)^{\mu(B^{(p)})}}} & \text{by Theorem \ref{thm.2}}
\\ &= 
\frac{
\lift_p \left(\intrat dX \,\, \twistrat_{\a_p} \Zratc(L) \right)}
{c_+^{\sigma_+( B^{(p)} ) } c_-^{\sigma_-( B^{(p)} ) }} & 
\\ &=
\frac{\int dX^{(p)} \,\,  \Zc(L^{(p)})}
{c_+^{\sigma_+( B^{(p)} ) } c_-^{\sigma_-( B^{(p)} ) }}
& \text{by Theorem \ref{thm.3}} 
\\ &=
\Z(\Sp_{(M,K)}). & \text{by $Z$'s definition}. 
\end{xalignat*}

\begin{theorem}
\lbl{thm.3}(Proof in Section \ref{sub.deflift})
For $\a_p=\nu^{-\frac{p-1}{p}}$ we have
$$
\lift_p \left(\intrat dX \,\, \twistrat_{\a_p}\Zratc(L) \right) 
= \int dX^{(p)} \,\,  \Zc(L^{(p)}).
$$
\end{theorem}

This reduces Theorem \ref{thm.1} to Theorems \ref{thm.2} and \ref{thm.3},
for a suitable $\lift_p$ map, and moreover, it shows that the presence of the 
\sig\ function in Theorem \ref{thm.1} is due to the normalization factors
$c_\pm$ of $\Zrat$.

The rest of the paper is devoted to the proof of Theorems \ref{thm.2} and 
\ref{thm.3} for a suitable residue map $\lift_p$. 

\section{Three views of knots}
\lbl{sec.comparison}

This section consists entirely of a classical topology  
view of knots and their abelian invariants such as \sig s,
Alexander polynomials and Blanchfield pairings. There is some overlap of
this section with \cite{GK1}; however for the benefit of the reader
we will try to present this section as self-contained as possible.

\subsection{The surgery and the Seifert surface view of knots}
\lbl{invariants}

In this section we discuss two views of knots $K$ in \ihs s $M$:
the {\em surgery view}, and the {\em Seifert surface} view.

We begin with the surgery view of knots. 
Given a surgery presentation $L$ for a pair $(M,K)$, let $W$ denote the
equivariant linking matrix of $L$, i.e., the linking matrix of a lift
$\ti L$ of $L$ to the universal cover $\widetilde{ST}$ of $ST$. 
It is not hard to see that $W$ is a Hermitian matrix, well-defined. Recall
the quotient $\B$ of the set of Hermitian matrices, from Section 
\ref{sec.review}. In \cite[Section 2]{GK2} it was shown that $W \in \B$ depends
only on the pair $(M,K)$ and not on the choice of a surgery presentation
of it. In addition, $W$ determines the {\em Blanchfield pairing} of $(M,K)$.
Thus, the natural map $\text{Knots}\to\text{BP}$ (where $\text{BP}$
stands for the {\em set of Blanchfield pairings})
factors through an (onto) map $\text{Knots}\to\B$. 

We now discuss the {\em Seifert surface view} of knots.
A more traditional way of looking at the set $\text{BP}$ of knots is
via Seifert surfaces and their associated Seifert matrices.
There is an onto map $\text{Knots}\to \Sei$,
where $\Sei$ is the set of matrices $A$ with integer entries
satisfying $\det(A-A')=1$, considered modulo an equivalence relation
called $S$-{\em equivalence}, \cite{Le}. It is known that the sets $\Sei$
and $\text{BP}$ are in 1-1 correspondence, see for example \cite{Le} and
\cite{Tr}.
Thus, we have a commutative diagram
$$
\divide\dgARROWLENGTH by2
\begin{diagram}
\node{\mathrm{Knots}}
\arrow{e,A}\arrow{s,A}
\node{\B}
\arrow{s,A} \\
\node{\Sei}
\arrow{e,t}{\sim}
\node{\text{BP}}
\end{diagram}
$$ 

It is well-known how to define abelian invariants of knots, such as 
the \sig\  and the {\em Alexander polynomial} $\Delta$, using Seifert surfaces.
Lesser known is a definition of these invariants using equivariant linking
martices, which we now give.

\begin{definition}
\lbl{def.W}
Let
$$
\d : \Herm(\Lz\to\BZ) \longto \Lz
$$
denote the (normalized) determinant given by $\d(W)=\det(W)\det(W(1))^{-1}$
(for all $W \in \Herm(\Lz\to\BZ)$) and let
$$
\vars: \Herm(\Lz\to\BZ) \longto \text{Maps}(S^1,\BZ)
$$
denote the function given by $\vars_z(W)=\s(W(z))-\s(W(1))$.
For a natural number $p$, let
$$
\vars_p: \Herm(\Lz\to\BZ) \longto \BZ
$$
be given by $\sum_{\w^p=1} \vars_{\w}(W)$.
\end{definition}

It is easy to see that $\d$ and $\vars$ descend to functions on $\B$. 
Furthermore, we have that
$$
\vars_p(W)=\s(W(T^{(p)}))-p\s(W(1)),
$$ where
$T^{(p)}$ is a $p$-cycle $p$ by $p$ matrix, given by example for $p=4$:
\begin{equation}
\lbl{eq.Tp}
T^{(4)} = \left[
\begin{array}{llll}
0 & 1 & 0 & 0 \\
0 & 0 & 1 & 0 \\
0 & 0 & 0 & 1 \\
1 & 0 & 0 & 0 
\end{array}
\right]
\end{equation}

\subsection{The clover view of knots}
\lbl{sub.Yview}

It seems hard to give an explicit algebraic map $\Sei \to\B$
although both sets may well be in 1-1 correspondence. Instead, we will
give a third view of knots, the {\em \clover\ view} of knots,
which enables us to prove Theorem \ref{thm.2}.

Consider a standard Seifert surface $\S$ of genus $g$ in $S^3$, which we 
think of as an embedded disk with pairs of bands attached in an alternating 
way along the disk:
$$
\psdraw{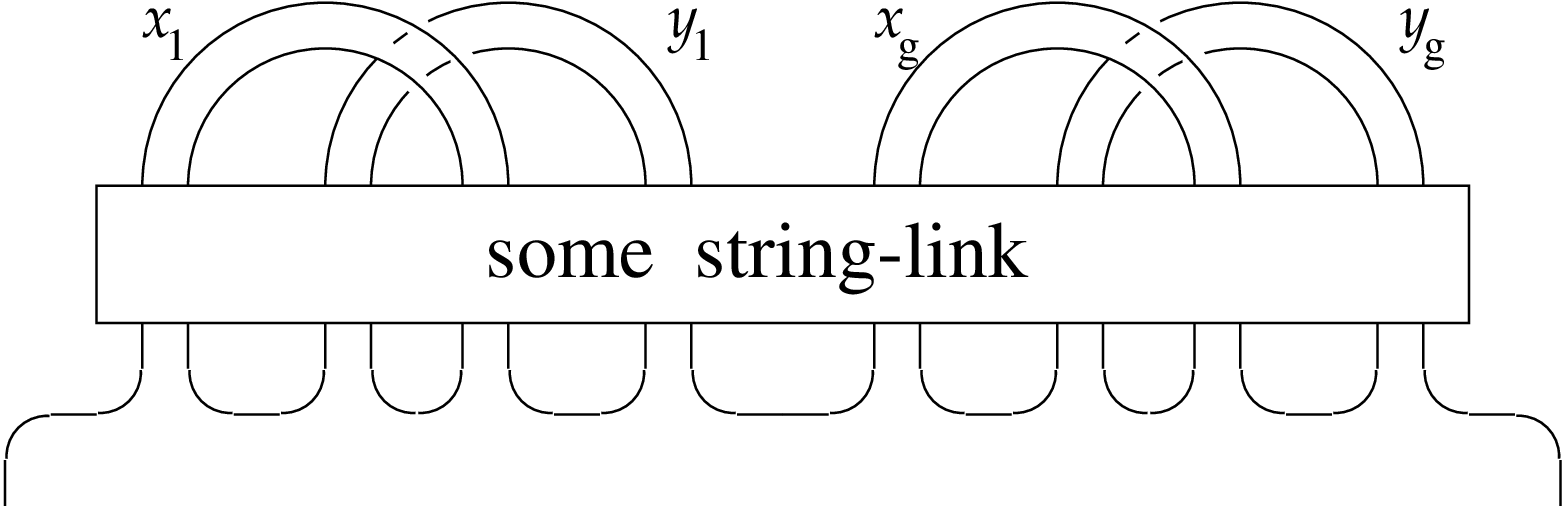}{2in}
$$ 
Consider an additional link $L'$ in $S^3\sminus \S$, such that its linking
matrix $C$ satisfies $\det(C)
=\pm 1$ and such that the linking number between the cores of the bands
and $L'$ vanishes. With respect to a suitable orientation 
of the 1-cycles corresponding to the cores of the bands, a Seifert matrix of 
$\S$ is given by
$$
A=
\left[
\begin{array}{ll}
L^{xx} & L^{xy} \\
L^{yx}-I & L^{yy}
\end{array}
\right], 
$$
where
$$
\left[
\begin{array}{ll}
L^{xx} & L^{xy} \\
L^{yx} & L^{yy}
\end{array}
\right] 
$$
is the linking matrix of the closure of the above string-link in the basis
$\{x_1,\dots,x_g,y_1,\dots,y_g\}$.
Let $(M,K)$ denote the pair obtained from $(S^3,\partial\S)$ after surgery on 
$L'$. With the notation  
$$
A \oplus B=
\left[
\begin{array}{ll}
A & 0 \\ 0 & B 
\end{array}
\right]
$$
we claim that 

\begin{theorem}
\lbl{thm.2g}
Given $(\S,L')$ as above, 
there exists a $2g$ component link $L$ in the complement of $L'$ such that:
 \newline
{\rm (a)}
$L \cup L' \subset ST $ is
a surgery presentation of $(M,K)$ in the sense of Section \ref{sub.surgery}.
\newline
{\rm (b)}
The equivariant linking matrix of $L \cup L'$ is represented by  
$W(t)\oplus C$ where
$$
W(t) = 
\left[
\begin{array}{llll}
L^{xx} & (1-t^{-1})L^{xy} - I \\
(1-t)L^{yx} - I& (1-t-t^{-1}+1)L^{yy} 
\end{array}
\right] .
$$
{\rm (c)}
Every pair $(M,K)$ comes from some $(\S,L')$ as above.
\end{theorem}

We will call such surgery presentations the {\em clover} view of knots.

\begin{proof}
(a) We will construct $L$ using the calculus of {\em clovers with two leaves}
introduced independently by Goussarov and Habiro \cite{Gu,Ha}; see also
\cite[Section 3]{GGP}. Clovers
with two leaves is a shorthand notation (on the left) for framed links  
shown on the right of the following figure:
$$
\psdraw{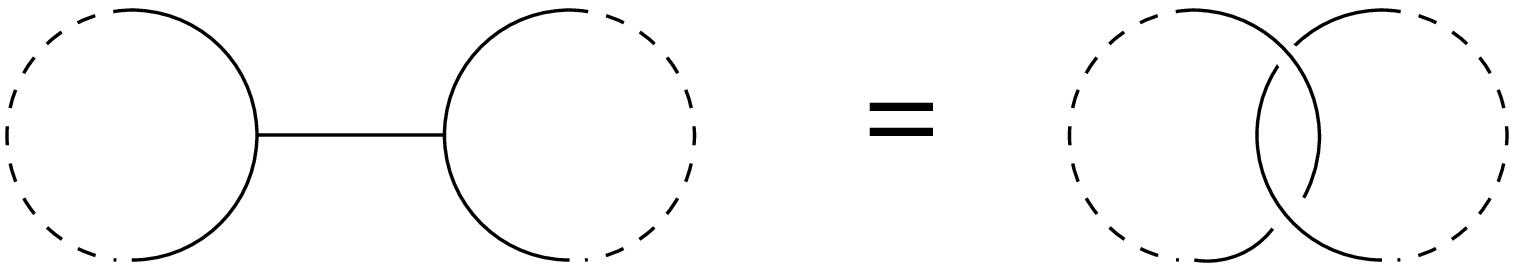}{1.5in}
$$
Since clovers can be thought of as framed links, surgery on clovers 
makes sense. Two clovers are equivalent (denoted by $\sim$ in the figures)
if after surgery, they represent the same 3-manifold. By calculus on clovers
(a variant of Kirby's calculus on framed links) we mean a set of moves 
that result to equivalent clovers. For an example of calculus on clovers,
we refer the reader to \cite{Gu,Ha} and also \cite[Sections 2,3]{GGP}.

In figures involving clovers, $L$ is constructed as follows:
$$
\psdraw{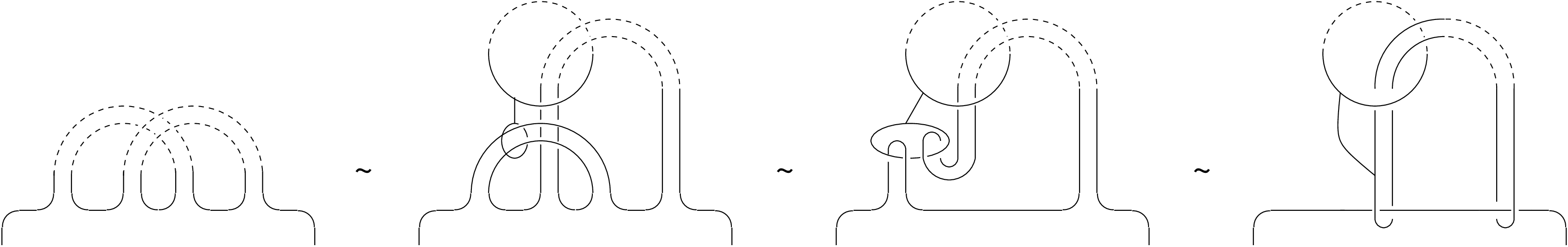}{5in}
$$
Notice that at the end of this construction, $L \cup L' \subset ST$
is a surgery presentation for $(M,K)$. \newline
(b) Using the discussion of \cite[Section 3.4]{Kr2}, it is easy to see that 
the equivariant linking matrix of (a based 
representative of) $L\cup L'$ is given as stated.
\newline
(c) Finally, we show that every pair $(M,K)$ arises this way. Indeed, choose
a Seifert surface $\S$ for $K$ in $M$ and a link $L' \subset M$ such that
$M_{L'}=S^3$. The link $L'$ may intersect $\S'$, and it may have
nontrivial linking number with the cores of the bands of $\S'$. However,
by a small isotopy of $L'$ in $M$ (which preserves the condition $M_{L'}=S^3$)
we can arrange that $L'$ be disjoint from $\S'$ and that its linking number 
with the cores of the bands vanishes. Viewed from $S^3$ (i.e., reversing
the surgery), this gives rise to $(\S,L')$ as needed.   
\end{proof}

The next theorem identifies the Alexander polynomial and the signature function
of a knot with the functions $\d$ and $\vars$ of Definition \ref{def.W}.

\begin{theorem}
\lbl{thm.Wd}
The maps composition of the maps $\d$ and $\vars$ with the natural map 
$\mathrm{Knots} \longto \B$ is given by the Alexander polynomial and the 
\sig\ function, respectively.
\end{theorem} 

\begin{proof}
There are several ways to prove this result, including an algebraic one,
which is a computation of appropriate Witt groups, and an analytic one,
which identifies the invariants with $U(1)$ $\rho$-invariants. None of
these proofs appear in the literature. We will give instead a proof using
the ideas already developed.

Fix a surgery presentation $L \cup L'$ for $(M,K)$, with equivariant linking 
matrix $W(t)\oplus C$ as in Theorem \ref{thm.2g}. 
Letting $P=\left[
\begin{array}{ll}
(1-t)I & 0 \\ 0 & I 
\end{array}
\right] \oplus I$, it follows that

\begin{eqnarray*}
P (W(t)\oplus C)P^\star & = & 
\left(\left[
\begin{array}{llll}
(1-t)I & 0 \\
0 & I      
\end{array}
\right] \oplus I \right)
\left(W(t) \oplus C \right)
\left(
\left[
\begin{array}{llll}
(1-t^{-1})I & 0  \\
0 & I   
\end{array}
\right]\oplus I \right) \\
& = &
\left[
\begin{array}{llll}
((1-t) + (1-t^{-1}))L^{xx}   &
((1-t) + (1-t^{-1}))L^{xy} - (1-t)I  \\
((1-t) + (1-t^{-1}))L^{yx} - (1-t^{-1})I &
((1-t) + (1-t^{-1}))L^{yy} 
\end{array}
\right] \oplus C\\
& = &
\left(
(1-t^{-1})A + (1-t)A'
\right) \oplus C.
\end{eqnarray*}
Taking signatures for any $t \in S^1$, $t \neq 1$, it follows that
\begin{eqnarray*}
\s(W(t))+\s(C) &=& \s(W(t)\oplus C)  \\
&=& \s(\left(
(1-t^{-1})A + (1-t)A'
\right) \oplus C)  \\
&=& \s(\left(
(1-t^{-1})A + (1-t)A'
\right))+\s(C) \\
&=& \s_t(M,K) + \s(C),
\end{eqnarray*}
where the last equality follows from the definition of the \sig,
see \cite[p. 289]{Ka} and \cite{Rf}. Thus, $\s(W(t))=\s_t(M,K)$. Since
$W(1)$ is a metabolic matrix, it follows that $\s(W(1))=0$, from which
it follows that $\vars(M,K)=\s(M,K)$. Taking determinants rather than 
signatures in the above discussion, it follows that $\d(M,K)=\D(M,K)$.
\end{proof}

\begin{proof}(of Theorem \ref{thm.2})
Fix a surgery presentation $L \cup L'$ for $(M,K)$, with equivariant linking 
matrix $W(t)\oplus C$ as in Theorem \ref{thm.2g}. Then the linking
matrix $B$ and $B^{(p)}$ of $L \cup L'$ and $L^{(p)} \cup L^{'(p)}$ are given 
by $W(1)\oplus C$ and $W(T^{(p)}) \oplus C \otimes I$ 
with an appropriate choice of basis. The result follows using Definition
\ref{def.W} and Theorem \ref{thm.Wd}.
\end{proof}

\begin{remark}
\lbl{rem.w}
An alternative proof of Theorem \ref{thm.2} can be obtained using the
$G$-signature theorem to the 4-manifold $N$ obtained by gluing two
4-manifolds $N_1,N_2$ with $\BZ_p$ actions along their common
boundary $\partial N_1=\partial N_2=\Sp_{(M,K)}$. Here $N_1$ is the
branched cover of $D^4$ branched along $D^2$ (obtained from adding
the handles of $L$ to $D^4$) and $N_2$ is a 4-manifold obtained from
a Seifert surface construction of $\Sp_{(M,K)}$.
\end{remark}

\begin{remark}
\lbl{rem.Y}
An alternative proof of Theorem \ref{thm.2g} can be obtained as follows. 
Start from a surgery presentation of $(M,K)$ in terms of clovers with 
three leaves, as was explained in \cite[Section 6.4]{GGP} and summarized
in the following figure:
$$ 
\psdraw{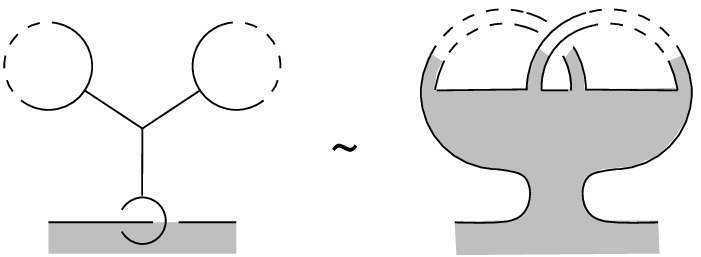}{2in} 
$$
Surgery on a clover with three leaves can be described in terms of surgery
on a six component link $L'''$. It was observed by the second author in 
\cite[Figure 3.1]{Kr3} that $L'''$ can be simplified via Kirby moves to a four 
component link $L''$. It is a pleasant exercise (left to the reader)
to further simplify $L''$ using Kirby moves to the two component link 
$L$ that appears in Theorem \ref{thm.2g}. 
\end{remark}

\begin{remark}
\lbl{rem.freedman}
Though we will not make use of this, we should mention that the clover
presentation $L$ of $(S^3,K)$ appears in work of M. Freedman
\cite[Lemma 1]{Fr}. Freedman starts with a knot of Arf invariant zero
together with a Seifert surface and constructs a spin 4-manifold $W_K$ 
with boundary $S^3_{K,0}$ (zero-surgery on $K$) by adding
suitable $1$-handles and $2$-handles in the $4$-ball.
The intersection form of $W_K$, as Freedman computes in \cite[Lemma 1]{Fr}
coincides with the equivariant linking matrix of $L$ of our
Theorem \ref{thm.2g}. This is not a coincidence, in fact the clover
view of knots,
interpreted in a 4-dimensional way as addition of $1$ and $2$ handles to
the $4$-ball, gives precisely Freedman's $4$-manifold.
\end{remark}

\section{Twisting}
\lbl{sec.twist}

In this section we define a notion of {\em twisting} $\twist_\a: 
\GA(\star_{X \cup k})\to \GA(\star_{X \cup k})$ and its rational cousin 
$\twistrat_\a:\GA(\star_X,\Lloc)\to
\GA(\star_X,\Lloc)$. Twisting (by elements of $\A(\star)$) is an operation on 
diagrams with beads which is analogous to the ``differential operator''
action of $\A(\star)$ on $\A(\star)$ defined in terms of gluing all legs
of the differential operator to some of the legs of a diagram.

A special case of twisting is the operation of wheeling on diagrams, studied
by \cite{A0,BLT,BL}. For a further discussion on the relation of
twisting and wheeling, see Section \ref{sec.wheel}. 

\subsection{Various kinds of diagrams}
\lbl{sub.various}

Manipulating the invariant $\Zrat$ involves calculations that take
values in vector spaces spanned by diagrams, modulo subspaces of relations.
The notation is as follows: given a ring $R$ with a distinguished group
of units $U$, and (possibly empty sets) $X,Y \cup T$, 
$\cD(\tline_X,\star_{Y \cup T},R,U)$ is the set of 
\begin{itemize}
\item
Uni-trivalent 
diagrams with skeleton $\tline_X$, with symmetric univalent vertices 
labeled by $Y \cup T$.
\item
The diagrams have oriented edges and skeleton and each edge is labeled
by an element of $R$, such that the edges that are part of the skeleton
are labeled only by $U$. Moreover, the product of the labels along each 
component
of the skeleton is $1$. Labels on edges or part of the skeleton will be 
called {\em beads}.
\end{itemize}

$\A(\tline_X,\star_Y,\ostar_T,R,U)$ is the quotient of the free vector space 
over $\BQ$
on $\cD(\tline_X,\star_{Y \cup T},R,U)$, modulo the relations of
\begin{itemize}
\item
$\AS$, $\IHX$, {\em multilinearity on the beads} shown in 
\cite[Figure 2]{GK2}.
\item
The {\em Holonomy Relation} shown in \cite[Figure 3]{GK2}.
\item
The {\em $T$-flavored basing relations} of \cite[Appendix D]{GK2}.
\end{itemize}

Empty sets will be omitted from the notation, and so will $U$, the selected
group of units of $R$. For example, $\A(\star_Y, R)$, $\A(R)$ and $\A(\phi)$
stands for $\A(\tline_\phi,\star_Y, \ostar_\phi,R,U)$, 
$\A(\tline_\phi,\star_\phi, \ostar_\phi,R,U)$ and $
\A(\tline_\phi,\star_\phi, \ostar_\phi,\BZ,1)$ respectively. 
Univalent vertices of diagrams will often be called {\em legs}. 
Diagrams will sometimes
be referred to as {\em graphs}. Special diagrams, called 
{\em struts}, labeled by $a,c$ with bead $b$ are drawn as follows
$$
\strutb{a}{c}{b}.
$$
oriented from bottom to top.

To further simplify notation, we will write $\A(\star), \A(\tline)$ and 
$\A(S^1)$ instead of $\A(\star_E), \A(\tline_E)$ and $\A(S^1_E)$
where $E$ is a set of one element.

A technical variant of the vector space $\A(\tline_X,\star_Y,\ostar_T,R,U)$ of
diagrams is the {\em set} $\GA(\tline_X,\star_Y,\ostar_T,R,U)$ 
which is the quotient of the set of group-like elements in
$\A(\tline_X,\star_Y,\ostar_T,R,U)$ (that is, exponential of a power series
of connected diagrams) modulo the group-like basing relation described in 
\cite[Section 3.3]{GK2}.

There is a natural map
$$
\GA(\tline_X,\star_Y,\ostar_T,R,U) \longrightarrow
\A(\tline_X,\star_Y,\ostar_T,R,U).
$$

Finally, $\GAz$ and $\A^0$ stand for $\B\times \GA$ and
$\B\times \GA$ respectively.

\subsection{A review of Wheels and Wheeling}
\lbl{sub.magic}
Twisting is closely related to the Wheels and Wheeling Conjectures introduced
in \cite{A0} and subsequently proven by \cite{BLT}. The Wheels and Wheeling 
Conjectures are a good tool to study structural properties of the Aarhus 
integral, as was explained in \cite{BL}. In our paper, they play a key role
in understanding twisting. 
In this section, we briefly review what Wheels and Wheeling is all about.

To warm up, recall that given an element $\a \in \A(\star)$ (such that
$\a$ does not contain a diagram one of whose components is a strut $\tline$)
we can turn it into an operator (i.e., linear map):
$$
\widehat{\a}: \, \A(\star)\to\A(\star)
$$
such that $\a$ acts on an element $x$ by gluing all legs of $\a$ to some of
the legs of $x$. It is easy to see that $\widehat{\a \sqcup \b}=
\widehat{\a} \circ \widehat{\b}$, which implies that if the constant term
of $\a$ is nonzero, then the operator $\a$ is invertible with
inverse $\widehat{\a}^{-1}=\widehat{\a^{-1}}$.

Of particular interest is the following element
$$
\Omega=\exp\left(\sum_{n=1}^\infty b_{2n} \wheel_{2n}\right) \in \A(\star)
$$
where $\wheel_{2n}$ is a wheel with $2n$ legs and 
$$
\sum_{n=1}^\infty b_{2n} x^{2n}=\frac{1}{2} \log \frac{\sinh x/2}{x/2} .
$$
The corresponding linear maps
$$
\homega^{-1}, \homega:\A(\star) \to\A(\star)
$$
are called respectively the {\em Wheeling} and the {\em Unwheeling}
maps and are denoted by $x \to x^{\wheel}$ and $x\to x^{\wheel^{-1}}$
respectively. Due to historical reasons dating
back to the days in Aarhus (where Wheeling was discovered) and also due to
Lie algebra reasons, wheeling was defined to be $\homega^{-1}$ and not 
$\homega$. 

Recall the symmetrization map $\chi: \A(\star)\to
\A(\tline)$ which sends an element $x \in \A(\star)$ to the average of 
the diagrams that arise by ordering the legs of $x$ on a line.
$\chi$ is a vector space isomorphism (with inverse $\s$) and can be 
used to transport the natural multiplication on $\A(\tline)$ (defined by 
joining two skeleton components of diagrams $\to \circ \to$ one next to the 
other to obtain a diagram on a skeleton component $\to$) to a multiplication 
on $\A(\star)$ which we denote by $\#$. There is
an additional multiplication $\sqcup$ on $\A(\star)$, defined using
the disjoint union of graphs. 

The {\em Wheeling Conjecture} states that the Unwheeling Isomorphism
$\homega: (\A(\ostar_k),\sqcup)\to (\A(\ostar_k),\#)$
interpolates the two multiplications on $\A(\ostar_k)$. Namely, that
for all $x,y \in \A(\ostar_k)$, we have
$$
\homega(x \sqcup y)=\homega(x) \, \# \, \homega(y).
$$

The {\em Wheels Conjecture} states that  
$$
\Z(S^3,\text{unknot})= \chi ( \Om) .
$$

The {\em Long Hopf Link Formula} states that
$$
Z\left(S^3,\psdraw{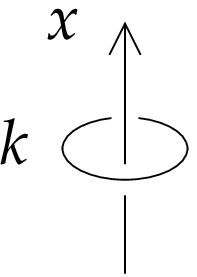}{0.25in}\right)
= 
\Om(k) \, \strutb{x}{}{e^k}  \in \A(\tline_x \ostar_k).
$$
Here and below, if $x \in \A(\star)$, then 
$x(h) \in \A(\star_h)$ denotes the diagram obtained from $x$ by replacing the 
color of the legs of $x$ by $h$.

It can be shown that the Wheels and Wheeling Conjectures are equivalent
to the Long Hopf Link Formula. In \cite{BLT} the Wheels and Wheeling 
Conjectures and the Long Hopf Link Formula were all proven. 
The identity $1+1=2$ (that is, doubling
the unknot component of the Long Hopf Link is a tangle isotopic to connecting 
sum twice the Long Hopf Link along the vertical strand), together with the
Long Hopf Link Formula imply the following {\em Magic Formula} 
\begin{equation}
\lbl{eq.magic}
\Om(k)\, \Om(h) \, \printname{3t}
	\setlength{\unitlength}{0.02\standardunitlength}
	\begin{array}{c}  \hspace{-1.7mm}
         	\raisebox{-8pt}{%\setlength{\unitlength}{0.00083333in}
\begingroup\makeatletter\ifx\SetFigFont\undefined%
\gdef\SetFigFont#1#2#3#4#5{%
  \reset@font\fontsize{#1}{#2pt}%
  \fontfamily{#3}\fontseries{#4}\fontshape{#5}%
  \selectfont}%
\fi\endgroup%
{\renewcommand{\dashlinestretch}{30}
\begin{picture}(711,1542)(0,-10)
%\put(75,312){\blacken\ellipse{120}{120}}
\put(15,312){$\bullet$}%\ellipse{120}{120}}
%\put(75,687){\blacken\ellipse{120}{120}}
\put(15,687){$\bullet$}%\ellipse{120}{120}}
\path(150,12)(150,1212)
\path(75,1137)(150,1212)(225,1137)
\put(450,837){\makebox(0,0)[lb]{$e^k$}}
\put(450,237){\makebox(0,0)[lb]{$e^h$}}
\put(0,1437){\makebox(0,0)[lb]{$x$}}
\end{picture}
}
 }
         	\hspace{-1.9mm}
	\end{array}
 =\Om(k+h) \strutb{x}{}{e^{k+h}}
\in \A(\tline_x, \ostar_{k,h})
\end{equation}

Before we end this section, we should mention that for $\a \in \A(\star)$,
the operator $\widehat{\a}$ can be defined for diagrams whose legs are 
colored by $X \cup \{k \}$ (abbreviated by $X \cup k$), where $k \not\in X$, 
by gluing all legs of $\a$ to some of the $k$-colored legs of a diagram. 
Furthermore, $\widehat{\a}$ preserves
$Y$-flavored basing relations for $Y \subset X \cup k$.
In addition, if $\a$ is group-like, then $\widehat{\a}$ sends group-like
elements to group-like elements. Note finally that $\A(\star)=\A(\ostar)$;
thus the operator $\widehat{\a}$ can be defined for $\a\in\A(\ostar)$.

\subsection{Twisting}
\lbl{sub.twist}

Throughout this section, $X$ denotes a (possibly empty) set disjoint from the 
two-element set $\{k,h\}$. Recall that given $x \in X$ and two diagrams 
$\a,\b \in \A(\star_X)$ with $k$ and $l$ $x$-colored legs respectively, 
the notation
$$
\la \a, \b \ra_{\{x\}} \in \A(\star_{X-\{x\}})
$$
means either zero (if $k \neq l$) or the sum of diagrams obtained by gluing
all $x$-colored legs of $\a$ with the $x$-colored legs of $\b$. This
definition can be extended to linear combination of diagrams, as a bilinear
symmetric operation, and can be further extended to an operation of gluing
$Y$-colored legs, for any $Y \subset X$.

\begin{remark}
We will often write 
$$
\la \a(y), \b(y) \ra_Y
$$
for the above operation, to emphasize the $Y$-colored legs of the diagrams.
{\em Warning:} In \cite{A0,GK2}, the authors used the alternative notation 
$\la \a(\pt y), \b(y)\ra_Y$ for the above operation.
\end{remark}

Given a diagram $s \in \A(\star_{X \cup k})$, 
the diagram $\phia(s) \in \A(\star_{X \cup k,h})$ denotes the sum of 
relabelings of legs of $s$ marked by $k$ by either $k$ or $h$. 

\begin{definition}
\lbl{def.twist}
For a group-like element $\a \in \A(\star)$, we 
define a map 
$$
\twist_\a: \A(\star_{X \cup k})\to \A(\star_{X \cup k})
$$
by 
$$
\twist_\a(s)=\la \phia(s) \Om(h)^{-1}, \a(h) \ra_h.
$$
\end{definition}

It is easy to see that $\twist_\a$ maps group-like elements to group-like
elements and maps $Y$-flavored basing relations to $Y$-flavored basing 
relations for $Y \subset X \cup k$; the latter follows from a ``sweeping
argument''.

The following lemma summarizes the elementary tricks about
the operators $\widehat{\a}$ that are very useful:

\begin{lemma}
\lbl{lem.trycks}
The operation $\la \cdot, \cdot \ra_X$ of gluing $X$-colored legs of diagrams
satisfies the following identities: 
$$
\la A(x), B(x) \sqcup C(x)\ra_X= \la \widehat{B} \, A( x), C(x)\ra_X= 
\la A(x+x'),B(x) \sqcup C(x')\ra_{X,X'}
$$
where $X'$ is a set in 1-1 correspondence with the set $X$.
\end{lemma}

In fact, twisting can be expressed in terms of the above action:

\begin{lemma}
\lbl{lem.twist} We have that:
\begin{eqnarray*}
\twist_\a &=& \widehat{\homega^{-1}(\a)} \\
\twist_\a \circ \twist_\b &=& \twist_{\a \# \b}
\end{eqnarray*}
\end{lemma}

\begin{proof}
Recall that $\widehat{\a}(y)=\la \a(h), y(k+h)\ra_h=\la y(k+h), \a(h)\ra_h$.
For the first part, we have: 

\begin{xalignat*}{2}
\twist_\a (x) &= \la x(k+h) \, \Om(h)^{-1}, \a(h) \ra_h & \\
&=
\la x(k+h), (\widehat{\Om}^{-1}(\a))(h)\ra_h & 
\text{by Lemma \ref{lem.trycks}}\\
&=
\widehat{\homega^{-1}(\a)}(x) & \text{by above discussion.}
\end{xalignat*}

For the second part, we have 
\begin{xalignat*}{2}
\twist_\a \circ \twist_\b &=
 \widehat{\homega^{-1}(\a)} \circ 
\widehat{\homega^{-1}(\b)}  & \\
&= \widehat{\homega^{-1}(\a) \sqcup \homega^{-1}(\b)} & \\
&= \widehat{\homega^{-1}(\a \# \b)} & \text{by Wheeling} \\
&= \twist_{\a \# \b}.  & %%\prtag{\Box}   
\eatline
\end{xalignat*}
\end{proof}

We now define a rational version 
$$\phim: \GA(\star_X, \Lloc)\to \GA(\star_{X \cup h},\Lloc)
$$ 
of the map $\phia$. The idea is that we substitute $te^h$ for $t$
(where $t$ and $h$ do not commute) and then replace $e^h$ by an exponential
of $h$-colored legs. This was explained in \cite[Section 3.1]{GK2} using
the notion of the {\em Cohn localization} of the free group in two generators.
We will not repeat the explanation of \cite{GK2} here, but instead use
the substitution map 
%$$
%\phim: \GA(\star_X, \Lloc)\to \GA(\star_{X \cup h},\Lloc)
%$$
freely. The reader may either refer to the above mentioned reference for a 
complete definition of the $\phim$ map, or may compromise with the following
property of the $\phim$ map:
$$
\phim \left( \strutb{}{}{p(t)/q(t)} \right)= \sum_{n=0}^\infty
\strutb{}{}{ p(te^h)/q(t) ((q(t)-q(te^h))/q(t))^n}.
$$
where $p,q \in \BQ[t^{\pm 1}]$ and $q(1)=\pm 1$.

\begin{definition}
\lbl{def.twistrat}
For a group-like element $\a \in \GA(\star)$, we 
define a map 
$$
\twistrat_\a: \GAz(\star_{X},\Lloc)\to \GAz(\star_{X},\Lloc)
$$
by 
$$
\twistrat_\a(M,s)=\left(M, \left\la \psi(M(t e^h)M(t)^{-1})
 % \exp\left( -\frac{1}{2} \tr\log(M(t e^h)M(t)^{-1})
%\right) 
\, \phim (s) 
, \a(h) \right\ra_h \right),
$$
where
$$
\psi(A)=\exp\left( -\frac{1}{2} \tr\log(A) \right).
$$
\end{definition}

\begin{remark}
\lbl{rem.substitution}
Here and below, we will be using the notation $\phi_{t\to e^k}(s)$
and $s(t \to e^k)$ to denote the substitution $t \to e^k$.
\end{remark}

The motivation for this rather strange definition comes from the proof of
Lemma \ref{lem.2twists} and Theorem \ref{prop.ztwist} below.

\begin{lemma}
\lbl{lem.t1}
$\twistrat_\a$ descends to a map:
$$
\GAz(\ostar_{X},\Lloc)\to \GAz(\ostar_{X},\Lloc)
$$
\end{lemma}

\begin{proof}
We need to show that the group-like basing relations are
preserved. With the notation and conventions of \cite[Section 3]{GK2},
there are two group-like basing relations $\bgp_1$ and $\bgp_2$ on
diagrams. It is easy to see that the $\bgp_1$ basing relation is preserved.
The $\bgp_2$ relation (denoted by $\stackrel{\bgp_2}\sim$) is generated in 
terms of a move of pushing $t$ on all legs (of some fixed color $x$) of a 
diagram. 
Given a diagram $s(x)$ with some $x$-colored
legs, let $s(xt)$ denote the result of pushing $t$ on every $x$-colored leg of 
$s(x)$. In order to show that the $\bgp_2$ relation is preserved, 
we need to show that $\twistrat_\a(M,s(xt)) \stackrel{\bgp}\sim 
\twistrat_\a(M,s(x))$. 

Ignoring the matrix part (i.e., setting $M$ the empty matrix), we can compute 
as follows: 
\begin{xalignat*}{2}
\twistrat_\a(s(xt)) &=
\la \phim(s)(\phim(xt)), \a(h)\ra_h & \\
&=
\la \phim(s)(xt e^{h'}), \a(h+h')\ra_{h,h'} & \text{by Lemma \ref{lem.trycks}}
\\
&\stackrel{\bgp_2}\sim 
\la  \phim(s)(x e^{h'}), \a(h+h')\ra_{h,h'} & \\
&\stackrel{\bgp_1}\sim 
\la  \phim(s)(x), \a(h+h')\ra_{h,h'} & \\
& = 
\la  \phim(s)(x), \a(h)\ra_{h,h'} & \\
& = 
\twistrat_\a(s(x))
\end{xalignat*}
The same calculation can be performed when we include the matrix part,
to conclude that $\twistrat_\a(M,s(xt)) \stackrel{\bgp}\sim 
\twistrat_\a(M,s(x))$.
\end{proof}

The next lemma about  $\twistrat$ should be compared with Lemma 
\ref{lem.twist} about  $\twist$: 

\begin{lemma}
\lbl{lem.twistrat}
We have
$$
\twistrat_\a \circ \twistrat_\b = \twistrat_{\a \# \b}
$$
\end{lemma}

\begin{proof}
Observe that
\begin{equation}
\lbl{eq.sigmachi}
\la e^h e^{h'}, \a(h) \sqcup \b(h')\ra_{h,h'}= \chi(\a) \# \chi(\b)
= \la e^h, \s(\chi(\a) \# \chi(\b)) \ra_h .
\end{equation}
In \cite[Section 3]{GK2}, it was shown that the ``determinant'' function $\psi$
is multiplicative, in the sense that (for suitable matrices $A,B$) we 
have:
\begin{equation}
\lbl{eq.chain}
\psi(A B)=\psi(A) \psi(B).
\end{equation}
Let us define $\pr: \GAz\to\GA$ to be the projection $(M,s)\to s$.
It suffices to show that $\pr \circ \twistrat_\a \circ \twistrat_\b = 
\pr \circ \twistrat_{\a \# \b}$. We compute this as follows:
\begin{xalignat*}{2}
\pr \circ \twistrat_{\a \# \b}(M,s) &=
\la \psi(M(t e^h)M(t)^{-1}) \, \phim (s) 
, (\a \# \b)(h) \ra_h  &  \\
&= 
\la \psi(M(t e^h e^{h'})M(t)^{-1}) \, \phimm (s) 
, \a(h) \sqcup \b(h') \ra_{h,h'} 
& \text{by \eqref{eq.sigmachi}} \\
&=
\la\la \psi(M(t e^h e^{h'})M(t e^{h'})^{-1})
\psi(M(t e^{h'})M(t)^{-1})  & \\
& \hspace{0.2cm} \phimm (s) 
, \a(h) \ra_h, \b(h') \ra_{h'} 
& \text{by \eqref{eq.chain}} \\
&=
\la \psi(M(t e^{h'})M(t)^{-1}) \phimp 
\la \psi(M(t e^h)M(t)^{-1}) \, \phim (s), \a(h)\ra_h , b(h')\ra_{h'} 
& \\
&=
\la \psi(M(t e^{h'})M(t)^{-1}) \phimp \circ
\pr \circ \twistrat_\a(M,s) ,\b(h')\ra_{h'} 
& \\
&=
\pr \circ \twistrat_\b ( \twistrat_\a (M,s)). & 
\end{xalignat*}
Since $\a \# \b = \b \# \a$, the result follows.
\end{proof}

Our next task is to relate the two notions $\twist,\twistrat$ of twisting.
In order to do so, recall the map
$$
\hair_k: \GA(\ostar_X,\Lloc)\to \GA(\ostar_{X \cup k})
$$
of \cite[Section 7.1]{GK2} defined by the substitution
$$
\strutb{}{}{t} \to \sum_{n=0}^\infty \frac{1}{n!} \printname{attacchn}
	\setlength{\unitlength}{0.03\standardunitlength}
	\begin{array}{c}  \hspace{-1.7mm}
         	\raisebox{-8pt}{\begingroup\makeatletter\ifx\SetFigFont\undefined%
\gdef\SetFigFont#1#2#3#4#5{%
  \reset@font\fontsize{#1}{#2pt}%
  \fontfamily{#3}\fontseries{#4}\fontshape{#5}%
  \selectfont}%
\fi\endgroup%
{\renewcommand{\dashlinestretch}{30}
\begin{picture}(1165,1239)(0,-10)
\path(12,12)(12,1212)
\path(42.000,1092.000)(12.000,1212.000)(-18.000,1092.000)
\path(12,912)(312,912)
\path(192.000,882.000)(312.000,912.000)(192.000,942.000)
\path(12,612)(312,612)
\path(192.000,582.000)(312.000,612.000)(192.000,642.000)
\path(12,312)(312,312)
\path(192.000,282.000)(312.000,312.000)(192.000,342.000)
\put(612,537){\makebox(0,0)[lb]{$n$ $h$-labeled legs}}
\end{picture}
} }
         	\hspace{-1.9mm}
	\end{array}
 
$$
and extended to a map
$$
\hairnu_k: \GAz(\ostar_X,\Lloc)\to \GA(\ostar_{X \cup k})
$$
by
$$
\hairnu_k(M,s)= \psi(M(e^k)M(1)^{-1})
%\exp_{\sqcup} 
%\left( -\frac{1}{2} \tr\log(M(e^k)M(1)^{-1}) \right)
\sqcup \hair_k(s) \sqcup \Om(k).
$$

Then, 

\begin{lemma}
\lbl{lem.2twists}
The following diagram commutes:
$$
%\divide\dgARROWLENGTH by2
\begin{diagram}
\node{\GA(\star)\times \GAz(\ostar_X,\Lloc)}
\arrow{e,t}{\twistrat}\arrow{s,l}{Id \times \hairnu_k}
\node{\GAz(\ostar_X,\Lloc)}
\arrow{s,r}{\hairnu_k} \\
\node{\GA(\star)\times \A(\ostar_{X \cup k})}
\arrow{e,t}{\twist}
\node{\A(\ostar_{X \cup k})}
\end{diagram}
$$ 
\end{lemma}

\begin{proof}
For $\a \in \GA(\star)$ and $(M,s) \in \GAz(\ostar_X,\Lloc)$, we have:
\begin{xalignat*}{2}
\twist \, (\hairnu_k(M,x)) &=
\la \hairnu_{k+h}(M,x) \Om^{-1}(h), \a(h)\ra_h & \\
&= \la \Om(k+h) \, \psi(M(e^{k+h})M(1)^{-1}) \,
x(t\to e^{k+h}) \, \Om^{-1}(h), \a(h) \ra_h & \\
&= \la \Om(k+h) \, \psi(M(e^{k+h})M(1)^{-1}) \,
x(t\to e^{k+h}), \widehat{\Om^{-1}}(\a)(h) \ra_h &
\text{by Lemma \ref{lem.trycks}} \\
&= \la \Om(k) \, \Om(h) \, 
\psi(M(e^ke^h)M(1)^{-1}) \, x(t\to e^k e^h), \widehat{\Om^{-1}}(\a)(h) \ra_h &
\text{by \eqref{eq.magic}} \\
&= \la \Om(k) \, \psi(M(e^ke^h)M(1)^{-1}) \, 
x(t\to e^k e^h), (\widehat{\Om} \, \widehat{\Om^{-1}})(\a)(h)\ra_h &
\text{by Lemma \ref{lem.trycks}} \\
&= \la \Om(k) \,  \psi(M(e^ke^h)M(1)^{-1}) \,
x(t\to e^k e^h), \a(h) \ra_h & \\
&= \Om(k) \, \psi(M(e^k)M(1)^{-1}) \, \phi_{t\to e^k}
\la \psi(M(te^h)M(t)^{-1}) \, x(t\to te^h), \a(h)\ra_h
& \text{by \eqref{eq.chain}} \\
&= \Om(k) \, \psi(M(e^k)M(1)^{-1}) \, \phi_{t\to e^k}
\pr \circ \twistrat_\a(M,x) 
& \text{by definition of $\twistrat$} \\
&= \hairnu_k (\twistrat_\a(s)). & %\prtag{\Box}
\end{xalignat*}
\end{proof}

The above lemma among other things explains the rather strange definition
of $\twistrat$.

\begin{corollary}
\lbl{cor.ZtoZ}
For all $\a \in \GA(\star)$ we have
$$
\hairnu \circ \twistrat_\a \circ \Zrat(M,K)= \twist_\a \circ \Z(M,K) 
\in \GA(\star).
$$
\end{corollary}

\begin{proof}
It follows from the above lemma, together with the fact that
$$
\hairnu \circ \Zrat(M,K)= \Z(M,K) \in \GA(\star),
$$
shown in \cite[Theorem 1.3]{GK2}.
\end{proof}

The next proposition states that $\twistrat$ intertwines (i.e., commutes with)
the integration map $\intrat$.

\begin{proposition}
\lbl{prop.ztwist}
For all $X' \subset X$ and $\a \in \GA(\star)$, 
the following diagram commutes:
$$
%\divide\dgARROWLENGTH by2
\begin{diagram}
\node{\GAz(\ostar_X,\Lloc)}
\arrow{e,t}{\intrat \,dX'}\arrow{s,l}{\twistrat_\a}
\node{\GAz(\ostar_{X-X'}\Lloc)}
\arrow{s,r}{\twistrat_\a} \\
\node{\GAz(\ostar_{X^{(p)}})}
\arrow{e,t}{\intrat \, dX'}
\node{\GAz(\ostar_{X-X'}\Lloc)}
\end{diagram}
$$ 
with the understanding that $\intrat$ is partially defined for $X'$-integrable
elements.
\end{proposition}

\begin{proof}
This is proven in \cite[Appendix E]{GK2} and repeated in Appendix 
\ref{sec.repeat}.
\end{proof}

\section{Lifting}
\lbl{sec.lifting}

%{\bf edit! from fileb.tex}
\subsection{The definition of the $\lift_p$ map}
\lbl{sub.deflift}

The goal of this section is to define the map $\lift_p$ and prove 
Theorem \ref{thm.3}.
We begin with a somewhat general situation. Consider a diagram $D$ with 
skeleton $\tline_X$, whose edges are labeled by elements of $\La$. 
For convenience, we express this by a diagram where there is a separate
bead for each $t^{\pm 1}$. $D$ consists of a solid part $\tline_X$ and a 
dashed part, that each have beads on them. 
The skeleton $X^{(p)}$ is defined by replacing each solid edge of $\tline_X$
by a parallel of $p$ solid edges. The skeleton $X^{(p)}$ has beads $t^{\pm 1}$
and the connected components of $X^{(p)}-(\text{beads})$ are labeled by 
$\BZ_p$ according to the figure shown below (for $p=4$)
$$
\psdraw{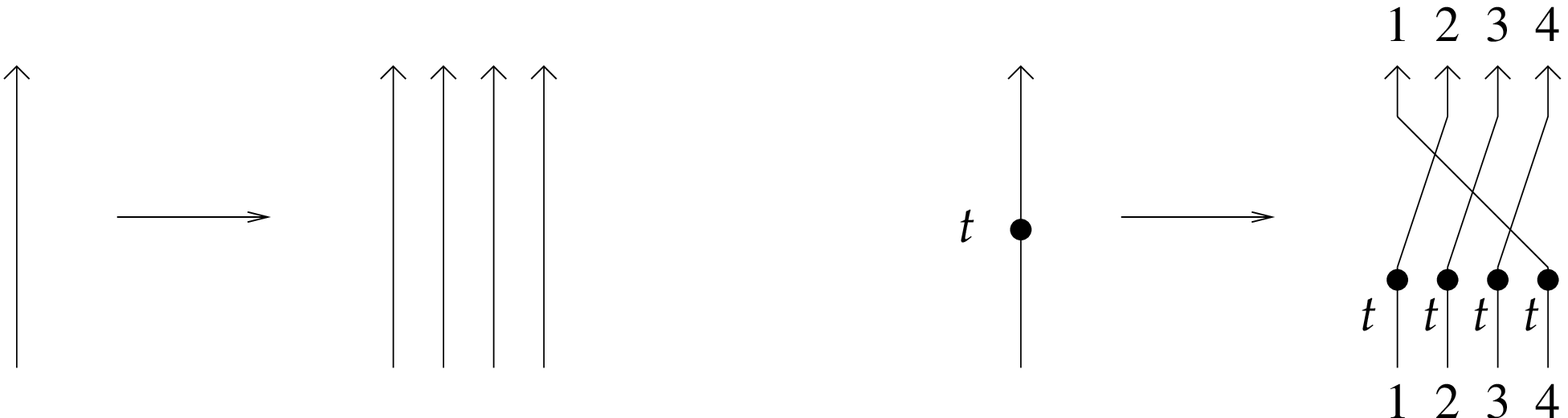}{3in}
$$
There is a projection map $\pi_p:X^{(p)}\to X$. A {\em lift} of a diagram
$D$ on $X$ is a diagram on $X^{(p)}$ whose dashed part is an isomorphic
copy of the dashed part of $D$, where the location on $X^{(p)}$ of each
univalent vertex maps under $\pi_p$ to the location of the corresponding 
univalent vertex on $X$. A $\BZ_p$-{\em labeling} of a diagram is 
an assignment of an element of $\BZ_p$ to each of the dashed or solid edges
that remain once we remove the beads of a diagram. A $\BZ_p$-labeling of
a diagram on $X^{(p)}$ is called $p$-{\em admissible}
if (after inserting the beads) it locally looks like
$$
\printname{p-coloring}
	\setlength{\unitlength}{0.03\standardunitlength}
	\begin{array}{c}  \hspace{-1.7mm}
         	\raisebox{-8pt}{\begingroup\makeatletter\ifx\SetFigFont\undefined%
\gdef\SetFigFont#1#2#3#4#5{%
  \reset@font\fontsize{#1}{#2pt}%
  \fontfamily{#3}\fontseries{#4}\fontshape{#5}%
  \selectfont}%
\fi\endgroup%
{\renewcommand{\dashlinestretch}{30}
\begin{picture}(6063,1539)(0,-10)
\put(5482,762){$\bullet$}%5562,762\ellipse{150}{150}}
\dashline{60.000}(762,1512)(762,12)
\dashline{60.000}(2412,1512)(3012,762)(3612,1512)
\dashline{60.000}(3012,12)(3012,762)
\dashline{60.000}(5562,1512)(5562,12)
\path(12,12)(1512,12)
\path(1392.000,-18.000)(1512.000,12.000)(1392.000,42.000)
\path(5487,1437)(5562,1512)(5637,1437)
\put(462,1062){\makebox(0,0)[lb]{$a$}}
\put(2337,1062){\makebox(0,0)[lb]{$a$}}
\put(3237,87){\makebox(0,0)[lb]{$a$}}
\put(5787,87){\makebox(0,0)[lb]{$a$}}
\put(5787,1287){\makebox(0,0)[lb]{$a+1$}}
\put(1212,237){\makebox(0,0)[lb]{$a$}}
\put(3537,1062){\makebox(0,0)[lb]{$a$}}
\put(237,237){\makebox(0,0)[lb]{$a$}}
\put(5787,687){\makebox(0,0)[lb]{$t$}}
\end{picture}
} }
         	\hspace{-1.9mm}
	\end{array}

$$ 
Now, we define $\lift_p(D)$ to be the sum of all diagrams on $X^{(p)}$ that
arise, when all the labels and beads are forgotten, from all $p$-admissible 
labelings of all lifts of $D$. As usual, the sum over the empty set is equal
to zero.

\begin{remark}
\lbl{rem.alt}
Here is an alternative description of $\lift_p(D,\a)$ for a labeling $\a$
of the edges of $D$ by monomials in $t$. Place a copy of $(D,\a)$ in $ST$ in
such a way that a bead $t$ corresponds to an edge going around the hole
of $ST$, as in \cite[Section 2.1]{Kr1}.
 Look at the $p$-fold cover $\pi_p: ST \to ST$, and consider
the preimage $\pi_p(D,\a) \subset ST \subset S^3$ as an abstract linear
combination of diagrams without beads. This linear
combination of diagrams equals to $\lift_p(D,\a)$.
\end{remark}

\begin{remark}
\lbl{rem.noskeleton}
Notice that in case $D$ has no skeleton, $b$ connected components, and all
the beads of its edges are $1$, then $\lift_p(D)=p^b D$.
\end{remark}

\begin{lemma}
\lbl{lem.liftp1}
The above construction gives a well-defined map
$$
\lift_p: \A(\tline_X,\La)\longrightarrow 
\A(\tline_{X^{(p)}})
$$
\end{lemma}

\begin{proof}
We need to show that the Holonomy Relations \cite[Figure 2]{GK2}
are preserved. There are two
possibilities: the case that all three edges in a Holonomy Relation are
dashed, and the case that two are part of the skeleton and the remaining
is dashed.

In the first case, the Holonomy Relation is preserved because there is an
obvious correspondence between lifts that admit an admissible labeling.

In the second case, the skeleton looks like (for $p=4$, with the convention
that $\bar t=t^{-1}$) 
$$
\psdraw{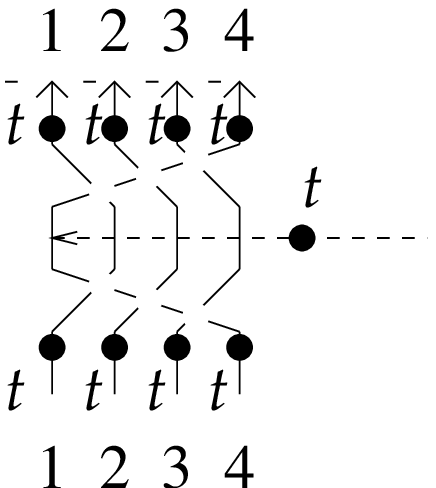}{1.0in}+
\psdraw{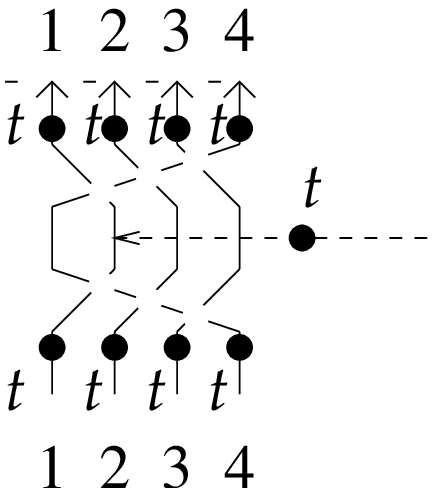}{1.0in}+
\psdraw{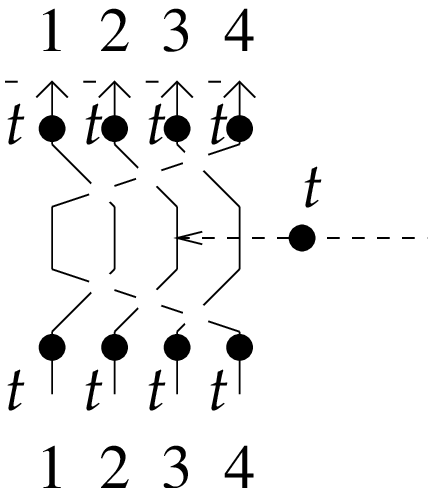}{1.0in}+
\psdraw{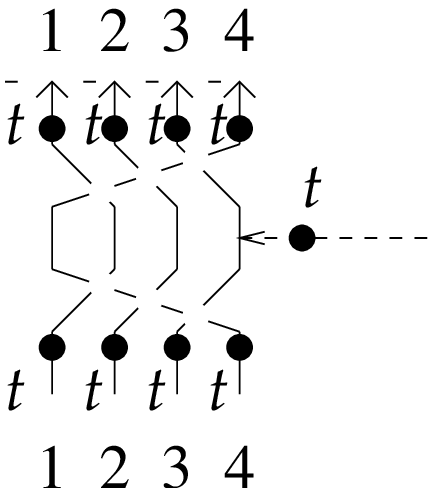}{1.0in}
% \,(\text{summed over lifts of}) \eepic{dline}{0.02}
$$
and again there is a correspondence between $p$-admissible labelings of
lifts of the two sides of the equation.
\end{proof}

There is a symmetrized version 
$$
\A(\star_X,\La)\longrightarrow \A(\star_{X^{(p)}})
$$
of the $\lift_p$ map, defined 
as follows: a lift of a diagram $D \in \A(\star_X,\La)$
is a diagram in $\A(\star_{X^{(p)}},\La)$ which consists of the same dashed
part as $D$, with each univalent vertex labeled by one of the $p$ copies of the
label of the univalent vertex of $D$ that it corresponds to. There is an
obvious notion of an admissible labeling of a diagram in 
$\A(\star_{X^{(p)}},\La)$, which is a labeling satisfying the conditions above,
and also 
$$
\printname{ddash}
	\setlength{\unitlength}{0.02\standardunitlength}
	\begin{array}{c}  \hspace{-1.7mm}
         	\raisebox{-8pt}{%\setlength{\unitlength}{0.00083333in}
\begingroup\makeatletter\ifx\SetFigFont\undefined%
\gdef\SetFigFont#1#2#3#4#5{%
  \reset@font\fontsize{#1}{#2pt}%
  \fontfamily{#3}\fontseries{#4}\fontshape{#5}%
  \selectfont}%
\fi\endgroup%
{\renewcommand{\dashlinestretch}{30}
\begin{picture}(974,1527)(0,-10)
\dashline{60.000}(150,1500)(150,300)
\put(375,900){\makebox(0,0)[lb]{$a \in \BZ_p$}}
\put(0,0){\makebox(0,0)[lb]{$x_i^{(a)}$}}
\end{picture}
}
 }
         	\hspace{-1.9mm}
	\end{array}
 \hspace{1cm}.
$$
Then, $\lift_p(D)$ is defined to be the sum of all diagrams on $X^{(p)}$ that
arise, when all the labels and beads are forgotten, from $p$-admissible 
labelings of lifts of $D$.

\begin{lemma}
\lbl{lem.liftp}
{\rm (a)} $\lift_p$ sends group-like elements to group-like elements and
induces maps that fit in the commutative diagram
$$
\divide\dgARROWLENGTH by2
\begin{diagram}
\node{\GA(\tline_X,\La)}
\arrow{e,t}{\s}\arrow{s,l}{\lift_p}
\node{\GA(\star_X,\La)}
\arrow{s,r}{\lift_p} \arrow{e}
\node{\GA(\ostar_X,\La)}
\arrow{s,r}{\lift_p} \\
\node{\GA(\tline_{X^{(p)}})}
\arrow{e,t}{\s}
\node{\GA(\star_{X^{(p)}})}
\arrow{e}
\node{\GA(\ostar_{X^{(p)}}).}
\end{diagram} 
$$
{\rm (b)} $\lift_p$ can be extended to a map
 $\GA(\ostar_X,\Lloc^{(p)})\to\GA(\ostar_{X^{(p)}})$,
where $\Lloc^{(p)}$ is the subring of $\Lloc$ that consists of all rational
functions whose denominators do not vanish at the complex $p$-th root of unity.
\end{lemma}

\begin{proof}
(a) Let us call an element of $\A(\tline_X,\La)$ {\em special} if the beads
of its skeleton equal to $1$. Using the Holonomy Relations, it follows that
$\A(\tline_X,\La)$ is spanned by special elements.

It is easy to see that $\lift_p$ maps group-like elements of
$\A(\star_X,\La)$ to group-like elements, and special group-like elements
in $\A(\tline_X,\La)$ to group-like elements in $\A(\tline_X,\La)$.  
Further, it is easy to show that the diagram
$$
\divide\dgARROWLENGTH by2
\begin{diagram}
\node{\A(\tline_X,\La)}
\arrow{e,t}{\s}\arrow{s,l}{\lift_p}
\node{\A(\star_X,\La)}
\arrow{s,r}{\lift_p} 
\arrow{s,r}{\lift_p} \\
\node{\A(\tline_{X^{(p)}})}
\arrow{e,t}{\s}
\node{\A(\star_{X^{(p)}})}
\end{diagram} 
$$
commutes when evaluated at special elements of $\A(\tline_X,\La)$.
From this, it follows that the left square diagram of the Lemma commutes.

For the right square, we need to show that the $X$-flavored {\em
basing relations} in $\GA(\star_X,\La)$
are mapped to $X^{(p)}$-flavored basing relations in $\GA(\star_{X^{(p)}})$.
There are two kinds of $X$-flavored basing relations, 
denoted by $\bgp_1$ and $\bgp_2$ in \cite[Section 3]{GK2}. First we consider
$\bgp_2$. Take two elements $s_1, s_2$ such that $s_1 
\stackrel{\bgp_2}\sim s_2$; we may assume that $s_2$
is obtained from pushing $t$ to each of the $x$-colored legs of $s_1$, for 
some $x \in X$. Corresponding to a diagram $D_1$ appearing in $s_1$, there
exists a diagram $D_2$ of $s_2$  obtained by pushing $t$ onto each of the 
$x$-colored legs of $D_1$. For example, 
$$
\printname{D1}
	\setlength{\unitlength}{0.03\standardunitlength}
	\begin{array}{c}  \hspace{-1.7mm}
         	\raisebox{-8pt}{\begingroup\makeatletter\ifx\SetFigFont\undefined%
\gdef\SetFigFont#1#2#3#4#5{%
  \reset@font\fontsize{#1}{#2pt}%
  \fontfamily{#3}\fontseries{#4}\fontshape{#5}%
  \selectfont}%
\fi\endgroup%
{\renewcommand{\dashlinestretch}{30}
\begin{picture}(2490,2277)(0,-10)
\put(117,1455){\arc{210}{1.5708}{3.1416}}
\put(117,2145){\arc{210}{3.1416}{4.7124}}
\put(1407,2145){\arc{210}{4.7124}{6.2832}}
\put(1407,1455){\arc{210}{0}{1.5708}}
\path(12,1455)(12,2145)
\path(117,2250)(1407,2250)
\path(1512,2145)(1512,1455)
\path(1407,1350)(117,1350)
\path(1512,1800)(2112,1800)
\path(312,450)(312,1350)
\path(342.000,1230.000)(312.000,1350.000)(282.000,1230.000)
\path(1212,450)(1212,1350)
\path(1242.000,1230.000)(1212.000,1350.000)(1182.000,1230.000)
\put(237,0){\makebox(0,0)[lb]{$x$}}
\put(1212,0){\makebox(0,0)[lb]{$x$}}
\put(2337,1725){\makebox(0,0)[lb]{$y$}}
\put(687,1725){\makebox(0,0)[lb]{$D_1$}}
\end{picture}
} }
         	\hspace{-1.9mm}
	\end{array}
 \hspace{1cm} \printname{D2}
	\setlength{\unitlength}{0.03\standardunitlength}
	\begin{array}{c}  \hspace{-1.7mm}
         	\raisebox{-8pt}{\begingroup\makeatletter\ifx\SetFigFont\undefined%
\gdef\SetFigFont#1#2#3#4#5{%
  \reset@font\fontsize{#1}{#2pt}%
  \fontfamily{#3}\fontseries{#4}\fontshape{#5}%
  \selectfont}%
\fi\endgroup%
{\renewcommand{\dashlinestretch}{30}
\begin{picture}(2553,2277)(0,-10)
%\put(375,825){\blacken\ellipse{120}{120}}
\put(275,825){$\bullet$}%\ellipse{120}{120}}
%\put(1275,825){\blacken\ellipse{120}{120}}
\put(1175,825){$\bullet$}%\ellipse{120}{120}}
\put(180,1455){\arc{210}{1.5708}{3.1416}}
\put(180,2145){\arc{210}{3.1416}{4.7124}}
\put(1470,2145){\arc{210}{4.7124}{6.2832}}
\put(1470,1455){\arc{210}{0}{1.5708}}
\path(75,1455)(75,2145)
\path(180,2250)(1470,2250)
\path(1575,2145)(1575,1455)
\path(1470,1350)(180,1350)
\path(1575,1800)(2175,1800)
\path(375,450)(375,1350)
\path(405.000,1230.000)(375.000,1350.000)(345.000,1230.000)
\path(1275,450)(1275,1350)
\path(1305.000,1230.000)(1275.000,1350.000)(1245.000,1230.000)
\put(300,0){\makebox(0,0)[lb]{$x$}}
\put(1275,0){\makebox(0,0)[lb]{$x$}}
\put(2400,1725){\makebox(0,0)[lb]{$y$}}
\put(750,1725){\makebox(0,0)[lb]{$D_2$}}
\put(1575,750){\makebox(0,0)[lb]{$t$}}
\put(0,750){\makebox(0,0)[lb]{$t$}}
\end{picture}
} }
         	\hspace{-1.9mm}
	\end{array}

$$
There is a 1-1 correspondence between admissible $p$-colorings of 
$\pi_p^{-1}(D_1)$ and those of $\pi_p^{-1}(D_2)$ (if we cyclically permute 
at the same time the labels $x^{(0)}, \dots, x^{(p-1)}$),  
shown as follows:
$$
\printname{D3}
	\setlength{\unitlength}{0.03\standardunitlength}
	\begin{array}{c}  \hspace{-1.7mm}
         	\raisebox{-8pt}{\begingroup\makeatletter\ifx\SetFigFont\undefined%
\gdef\SetFigFont#1#2#3#4#5{%
  \reset@font\fontsize{#1}{#2pt}%
  \fontfamily{#3}\fontseries{#4}\fontshape{#5}%
  \selectfont}%
\fi\endgroup%
{\renewcommand{\dashlinestretch}{30}
\begin{picture}(2492,2277)(0,-10)
\put(180,1455){\arc{210}{1.5708}{3.1416}}
\put(180,2145){\arc{210}{3.1416}{4.7124}}
\put(1470,2145){\arc{210}{4.7124}{6.2832}}
\put(1470,1455){\arc{210}{0}{1.5708}}
\path(75,1455)(75,2145)
\path(180,2250)(1470,2250)
\path(1575,2145)(1575,1455)
\path(1470,1350)(180,1350)
\path(1575,1800)(2175,1800)
\path(375,450)(375,1350)
\path(405.000,1230.000)(375.000,1350.000)(345.000,1230.000)
\path(1275,450)(1275,1350)
\path(1305.000,1230.000)(1275.000,1350.000)(1245.000,1230.000)
\put(300,0){\makebox(0,0)[lb]{$x^{(r)}$}}
\put(1275,0){\makebox(0,0)[lb]{$x^{(s)}$}}
\put(2400,1725){\makebox(0,0)[lb]{$y^{(k)}$}}
\put(750,1725){\makebox(0,0)[lb]{$D_1$}}
\put(0,750){\makebox(0,0)[lb]{$r$}}
\put(1500,750){\makebox(0,0)[lb]{$s$}}
\end{picture}
} }
         	\hspace{-1.9mm}
	\end{array}
 \longleftrightarrow \hspace{0.5cm} \printname{D4}
	\setlength{\unitlength}{0.03\standardunitlength}
	\begin{array}{c}  \hspace{-1.7mm}
         	\raisebox{-8pt}{\begingroup\makeatletter\ifx\SetFigFont\undefined%
\gdef\SetFigFont#1#2#3#4#5{%
  \reset@font\fontsize{#1}{#2pt}%
  \fontfamily{#3}\fontseries{#4}\fontshape{#5}%
  \selectfont}%
\fi\endgroup%
{\renewcommand{\dashlinestretch}{30}
\begin{picture}(2492,2277)(0,-10)
%\put(375,825){\blacken\ellipse{120}{120}}
\put(275,825){$\bullet$}%\ellipse{120}{120}}
%\put(1275,825){\blacken\ellipse{120}{120}}
\put(1175,825){$\bullet$}%\ellipse{120}{120}}
\put(180,1455){\arc{210}{1.5708}{3.1416}}
\put(180,2145){\arc{210}{3.1416}{4.7124}}
\put(1470,2145){\arc{210}{4.7124}{6.2832}}
\put(1470,1455){\arc{210}{0}{1.5708}}
\path(75,1455)(75,2145)
\path(180,2250)(1470,2250)
\path(1575,2145)(1575,1455)
\path(1470,1350)(180,1350)
\path(1575,1800)(2175,1800)
\path(375,450)(375,1350)
\path(405.000,1230.000)(375.000,1350.000)(345.000,1230.000)
\path(1275,450)(1275,1350)
\path(1305.000,1230.000)(1275.000,1350.000)(1245.000,1230.000)
\put(300,0){\makebox(0,0)[lb]{$x^{(r-1)}$}}
\put(1275,0){\makebox(0,0)[lb]{$x^{(s-1)}$}}
\put(2400,1725){\makebox(0,0)[lb]{$y^{(k)}$}}
\put(750,1725){\makebox(0,0)[lb]{$D_2$}}
\put(1675,750){\makebox(0,0)[lb]{$t$}}
\put(-200,750){\makebox(0,0)[lb]{$t$}}
\put(75,975){\makebox(0,0)[lb]{$r$}}
\put(-600,450){\makebox(0,0)[lb]{$r-1$}}
\put(1500,975){\makebox(0,0)[lb]{$s$}}
\put(1500,450){\makebox(0,0)[lb]{$s-1$}}
\end{picture}
} }
         	\hspace{-1.9mm}
	\end{array}
.
$$
Applying $\bgp_2$ basing relations, the two results agree.
In other words, $\lift_p(D) \stackrel{\bgp_2}\sim \lift_p(D')$.
 
Now, consider the case of $\bgp_1$, (in the formulation of 
\cite[Section 3]{GK2}). Given $s_1 \stackrel{\bgp_1}\sim s_2$,
there exists an element $s \in \GA(\star_{X \cup \pt h},
\La)$ with some legs labeled by $\pt h$, such that
\begin{eqnarray*}
s_1 &=& \conh (s) \\
s_2 &=& \conh (s(x \to x e^h))
\end{eqnarray*}
for some $x \in X$, where $\conh$ is the operation that contracts all $\pt h$
legs of a diagram to all $h$ legs of it. Now observe that 
\begin{eqnarray*}
\lift_p(s_2) &=& \lift_p( \conh(s(x \to xe^h))) \\
&=&
\mathrm{con}_{\{ h^{(0)}, \dots, h^{(p-1)}\}} \circ \lift_p (s(x \to xe^h)) \\
&=& 
\mathrm{con}_{\{ h^{(0)}, \dots, h^{(p-1)}\}} \circ \lift_p
s(x^{(0)} \to x^{(0)}e^{h^{(0)}}, \dots,x^{(p-1)} \to x^{(p-1)}
e^{h^{(p-1)}} ) \\
&\stackrel{\bgp_1}\sim &
\lift_p (s (h \to 0)) \\
&=& 
\lift_p(s_1).
\end{eqnarray*}
(b) Notice first that $\lift_p$ can be defined when beads
are labeled by elements of $\BC[t]/(t^p-1)$.
There is an isomorphism $ \Lloc^{(p)}/(t^p-1) \cong \BC[t]/(t^p-1)$ over $\BC$
which gives rise (after composition with the projection $ \Lloc^{(p)}\to
\Lloc^{(p)}/(t^p-1)$) to a map 
\begin{equation}
\lbl{eq.chp}
\ch_p:\Lloc^{(p)}\to \BC[t]/(t^p-1).
\end{equation}
Using this map, we can define $\lift_p$
as before and check that the relations are preserved.
\end{proof}

\begin{remark}
\lbl{rem.lft}
$\lift_p$ can also be extended to a map
$$
\lift_p: \GAz(\ostar_X,\Lloc^{(p)}) \longrightarrow\GA(\ostar_{X^{(p)}})
$$
by forgetting the matrix part, i.e., by $\lift_p(M,s)=\lift_p(s)$.
\end{remark}

Let $L$ be a surgery presentation of a pair $(M,K)$ as in Section 
\ref{sub.surgery} and let $L^{(p)}$ be the lift of $L$
to the $p$-fold cover of the solid torus, regarded as a link in $S^3$.
The following proposition is a key point.

\begin{proposition}
\lbl{prop.fundlem}
With the above notation, we have:
$$
\Zc(L^{(p)}) = \lift_p \circ \twistrat_{\a_p} \circ \Zratc(L).
$$
\end{proposition}

{\em Proof.}
We begin by recalling first how $\Zrat(L)$ is defined, following 
\cite[Section 4]{GK2}. The definition is given by representing $L$
in terms of objects called {\em sliced crossed links} in a solid torus. 
Sliced crossed links are planar tangles of a specific shape that can be 
obtained from a generic height function of a link $L$ in a standard solid torus
$ST$. Each component of their corresponding link in $ST$ is marked by
a cross ($\times$). 
Given a null homotopic link $L$ in $ST$, choose a sliced crossed link
representative $(T_0,T_1,T_2)$ where $T_0$ consists of local minima, $T_2$
consists of local maxima and $T_1$, thought of as a tangle in $I \times I$,
equals to $I_w \sqcup \gamma$. 
Here $w$, the {\em gluing site}, is a sequence in $\uparrow$ and $\downarrow$,
and $\bar{w}$ is the reverse sequence (where the reverse of $\uparrow\uparrow
\downarrow$ is $\downarrow\uparrow\uparrow$).

For example, for $w=\downarrow\uparrow$, we may have the following presentation
of a knot in $ST$
$$
\printname{figure8}
	\setlength{\unitlength}{0.025\standardunitlength}
	\begin{array}{c}  \hspace{-1.7mm}
         	\raisebox{-8pt}{\begingroup\makeatletter\ifx\SetFigFont\undefined%
\gdef\SetFigFont#1#2#3#4#5{%
  \reset@font\fontsize{#1}{#2pt}%
  \fontfamily{#3}\fontseries{#4}\fontshape{#5}%
  \selectfont}%
\fi\endgroup%
{\renewcommand{\dashlinestretch}{30}
\begin{picture}(2412,4389)(0,-10)
\path(675,3237)(825,3087)
\path(675,3087)(825,3237)
\put(1500.000,3437.000){\arc{650.000}{2.7468}{6.6780}}
\put(1425.000,3387.000){\arc{1358.308}{3.0309}{6.3938}}
\path(600,4062)(600,4362)(2400,4362)(2400,4062)
\path(600,4062)(600,3312)(2400,3312)(2400,4137)
\put(2288.000,1699.000){\arc{150.000}{6.2832}{9.4248}}
\put(2306.246,1717.247){\arc{119.230}{3.4793}{6.5943}}
\path(2213,1699)(2213,1849)
\put(2288.000,1399.000){\arc{150.000}{6.2832}{9.4248}}
\put(2306.246,1417.247){\arc{119.230}{3.4793}{6.5943}}
\path(2213,1399)(2213,1549)
\put(2288.000,1962.000){\arc{150.000}{6.2832}{9.4248}}
\put(2306.500,1980.500){\arc{118.903}{3.4580}{6.5996}}
\path(2213,1962)(2213,1963)(2212,1967)
	(2211,1977)(2209,1989)(2208,2003)
	(2208,2016)(2209,2027)(2213,2037)
	(2218,2045)(2224,2053)(2233,2063)
	(2245,2073)(2257,2084)(2269,2095)
	(2278,2104)(2285,2109)(2288,2112)
\path(2213,1324)(2213,1323)(2212,1318)
	(2211,1308)(2209,1295)(2208,1280)
	(2208,1267)(2209,1257)(2213,1249)
	(2219,1242)(2229,1236)(2242,1230)
	(2257,1224)(2271,1218)(2282,1214)
	(2287,1212)(2288,1212)
\put(1500.000,2112.000){\arc{600.000}{6.2832}{9.4248}}
\path(750,1212)(750,2112)
\path(1950,2112)(1950,1212)
\path(1080.000,1992.000)(1050.000,2112.000)(1020.000,1992.000)
\path(1050,2112)(1050,1212)
\path(600,2112)(600,1212)(2400,1212)(2400,2112)
\put(1500.000,992.357){\arc{910.714}{6.1296}{9.5783}}
\put(1500.000,1062.000){\arc{1500.000}{6.2832}{9.4248}}
\path(600,1062)(2400,1062)(2400,12)
	(600,12)(600,1062)
\put(1500.000,2262.000){\arc{600.000}{3.1416}{6.2832}}
\put(2100.000,2712.000){\arc{300.000}{3.1416}{6.2832}}
\put(1950.000,2787.000){\arc{300.000}{6.2832}{9.4248}}
\path(750,2262)(750,3162)
\path(1950,2262)(1950,2562)
\path(2250,2712)(2250,2262)
\path(1800,2787)(1800,3162)
\path(2100,2937)(2100,3162)
\path(600,3162)(2400,3162)
\path(600,3162)(600,2262)
\path(2400,3162)(2400,2262)
\path(600,2262)(1200,2262)
\path(1800,2262)(2400,2262)
\path(1800,2112)(2400,2112)
\path(600,2112)(1200,2112)
\path(1050,2262)(1050,2265)(1049,2271)
	(1048,2283)(1047,2299)(1046,2320)
	(1044,2345)(1043,2372)(1041,2400)
	(1040,2428)(1040,2455)(1040,2479)
	(1041,2502)(1043,2523)(1046,2543)
	(1050,2562)(1055,2579)(1060,2595)
	(1067,2612)(1075,2629)(1084,2645)
	(1093,2662)(1103,2679)(1114,2695)
	(1125,2712)(1136,2729)(1147,2745)
	(1157,2762)(1166,2779)(1175,2795)
	(1183,2812)(1190,2829)(1195,2845)
	(1200,2862)(1204,2881)(1207,2901)
	(1209,2922)(1210,2945)(1210,2969)
	(1210,2996)(1209,3024)(1208,3052)
	(1206,3079)(1204,3104)(1203,3125)
	(1202,3141)(1201,3153)(1200,3159)(1200,3162)
\put(0,3837){\makebox(0,0)[lb]{$T_2$}}
\put(0,2262){\makebox(0,0)[lb]{$T_1$}}
\put(0,537){\makebox(0,0)[lb]{$T_0$}}
\end{picture}
} }
         	\hspace{-1.9mm}
	\end{array}
 \hspace{1cm} \text{where} \hspace{1cm}
\gamma=\printname{gamma}
	\setlength{\unitlength}{0.02\standardunitlength}
	\begin{array}{c}  \hspace{-1.7mm}
         	\raisebox{-8pt}{\begingroup\makeatletter\ifx\SetFigFont\undefined%
\gdef\SetFigFont#1#2#3#4#5{%
  \reset@font\fontsize{#1}{#2pt}%
  \fontfamily{#3}\fontseries{#4}\fontshape{#5}%
  \selectfont}%
\fi\endgroup%
{\renewcommand{\dashlinestretch}{30}
\begin{picture}(924,1764)(0,-10)
\put(801.500,464.794){\arc{149.618}{0.0909}{3.0506}}
\put(817.913,470.375){\arc{118.781}{3.5324}{6.4931}}
\path(727,458)(727,596)
\put(801.500,189.794){\arc{149.618}{0.0909}{3.0506}}
\put(817.913,195.375){\arc{118.781}{3.5324}{6.4931}}
\path(727,183)(727,321)
\put(801.500,705.794){\arc{149.618}{0.0909}{3.0506}}
\put(818.206,711.697){\arc{118.345}{3.5099}{6.4994}}
\path(727,699)(727,700)(726,704)
	(725,712)(723,724)(722,737)
	(722,749)(724,759)(727,768)
	(732,775)(738,783)(747,792)
	(758,801)(770,812)(781,821)
	(791,829)(797,834)(800,837)
\path(727,115)(727,114)(726,110)
	(725,101)(723,88)(722,75)
	(722,63)(724,53)(727,46)
	(733,40)(742,34)(755,29)
	(770,23)(784,18)(794,14)
	(799,12)(800,12)
\put(612.000,1287.000){\arc{300.000}{3.1416}{6.2832}}
\put(462.000,1362.000){\arc{300.000}{6.2832}{9.4248}}
\path(462,837)(462,1137)
\path(762,1287)(762,837)
\path(312,1362)(312,1737)
\path(612,1512)(612,1737)
\path(12,1737)(912,1737)
\path(912,1737)(912,837)
\path(463,837)(463,12)
\path(12,537)(12,12)(912,12)(912,837)
\path(12,537)(12,1737)
\end{picture}
} }
         	\hspace{-1.9mm}
	\end{array}

$$
(and where the sliced crossed link is a tangle in an annulus). 
For typographical reasons, we will often say that $(T_0,T_1,T_2)$ is
the {\em closure} of the tangle $\gamma$. 

Consider a representation of a null homotopic link $L$ in $ST$
by $(T_0,T_1,T_2)$ as above. Recall that the fractional powers of $\nu$ in 
the algebra $(A(\ostar), \#)$ are defined as follows: for integers $n,m$,
$\nu^{n/m} \in A(\ostar)$ is the unique element whose constant term is $1$
that satisfies $(\nu^{n/m})^m=\nu^n$. 
 
Then, $\Zrat(L)$ is defined as the element
of $\GA(\ostar_X,\La)$ obtained by composition of
$$ (Z(T_0), I_{\bar w}(1) \otimes \D_{w}(\nu^{1/2}), 
I_{\bar w}(1) \otimes I_w(t), I_{\bar w}(1) \otimes 
Z(\gamma), I_{\bar w}(1) \otimes \D_w(\nu^{1/2}), Z(T_2))
$$
where $I_w(a)$ means a skeleton component that consists of solid arcs with
orientations according to the arrows in $w$, with $a$ (resp. $\bar{a}$)
placed on each $\uparrow$ (resp. $\downarrow$), and $\D_w$ is the 
commultiplication obtained by replacing a solid segment $\tline$ by
a $w$-parallel of it. After cutting the sliced crossed link at the crosses
($\times$), we consider the resulting composition of diagrams as an element of 
$\GA(\ostar_X,\La)$. We claim that
 
\begin{lemma}
\lbl{lem.tZ}
$\twistrat_\a \circ \Zrat(L) \in \GA(\ostar_X,\La)$ equals to the element
obtained by composition of
$$ (Z(T_0), I_{\bar w}(1) \otimes \D_{w}(\nu^{1/2}), I_{\bar w}(1) 
\otimes \D_w(\a), I_{\bar w}(1) \otimes I_w(t), I_{\bar w}(1) \otimes 
Z(\gamma), I_{\bar w}(1) \otimes \D_w(\nu^{1/2}), Z(T_2))
$$
\end{lemma}

\begin{proof}
This follows easily from the definition of the $\twistrat_\a$ using the
fact that the beads of the diagrams in $\Zratc(L)$ appear only at the gluing
site.
\end{proof}

In short, we will say that $\twistrat_\a \circ \Zrat(L)$ is obtained by the 
{\em closure} of the following sequence 
$$
(\D_{w}(\nu^{1/2}), \D_w(\a), I_w(t),  Z(\gamma), 
\D_w(\nu^{1/2})),
$$
which we will draw schematically as follows:
$$
\printname{tangle2a}
	\setlength{\unitlength}{0.02\standardunitlength}
	\begin{array}{c}  \hspace{-1.7mm}
         	\raisebox{-8pt}{\begingroup\makeatletter\ifx\SetFigFont\undefined%
\gdef\SetFigFont#1#2#3#4#5{%
  \reset@font\fontsize{#1}{#2pt}%
  \fontfamily{#3}\fontseries{#4}\fontshape{#5}%
  \selectfont}%
\fi\endgroup%
{\renewcommand{\dashlinestretch}{30}
\begin{picture}(2424,6039)(0,-10)
\put(687,3012){$\bullet$}%\ellipse{150}{150}}
\put(1587,3012){$\bullet$}%\ellipse{150}{150}}
\path(762,6012)(762,5712)
\path(1662,6012)(1662,5712)
\path(762,4512)(762,4812)
\path(1662,4512)(1662,4812)
\path(762,3612)(762,3012)
\path(1662,3012)(1662,2412)
\path(762,1512)(762,1212)
\path(1662,1512)(1662,1212)
\path(1662,312)(1662,12)
\path(762,312)(762,12)
\path(1662,3612)(1662,3012)
\path(1632.000,3132.000)(1662.000,3012.000)(1692.000,3132.000)
\path(792.000,2892.000)(762.000,3012.000)(732.000,2892.000)
\path(762,3012)(762,2412)
\path(12,5712)(2412,5712)(2412,4812)
	(12,4812)(12,5712)
\path(12,4512)(2412,4512)(2412,3612)
	(12,3612)(12,4512)
\path(12,2412)(2412,2412)(2412,1512)
	(12,1512)(12,2412)
\path(12,1212)(2412,1212)(2412,312)
	(12,312)(12,1212)
\put(762,3762){\makebox(0,0)[lb]{$\Z(\gamma)$}}
\put(1012,2937){\makebox(0,0)[lb]{$\bar t$}}
\put(2012,2937){\makebox(0,0)[lb]{$ t$}}
\put(537,4962){\makebox(0,0)[lb]{$\D(\nu^{1/2})$}}
\put(387,1662){\makebox(0,0)[lb]{$\D(\a)$}}
\put(537,462){\makebox(0,0)[lb]{$\D(\nu^{1/2})$}}
\end{picture}
} }
         	\hspace{-1.9mm}
	\end{array}

$$
Going back to the proof of Proposition \ref{prop.fundlem}, using
$\a_p=\nu^{-(p-1)/p}$, and the group-like basing relations on 
$\GA(\ostar_X,\La)$, it follows that we can slide and cancel the powers of 
$\nu$. Thus the closure of the above sequence for $\a=\a_p$,
equals to the following sequence:
$$
\printname{tangle2b}
	\setlength{\unitlength}{0.02\standardunitlength}
	\begin{array}{c}  \hspace{-1.7mm}
         	\raisebox{-8pt}{\begingroup\makeatletter\ifx\SetFigFont\undefined%
\gdef\SetFigFont#1#2#3#4#5{%
  \reset@font\fontsize{#1}{#2pt}%
  \fontfamily{#3}\fontseries{#4}\fontshape{#5}%
  \selectfont}%
\fi\endgroup%
{\renewcommand{\dashlinestretch}{30}
\begin{picture}(2424,3639)(0,-10)
\put(687,1812){$\bullet$}%\ellipse{150}{150}}
\put(1587,1812){$\bullet$}%\ellipse{150}{150}}
\path(762,312)(762,12)
\path(1662,312)(1662,12)
\path(762,3312)(762,3612)
\path(1662,3312)(1662,3612)
\path(12,3312)(2412,3312)(2412,2412)
	(12,2412)(12,3312)
\path(12,1212)(2412,1212)(2412,312)
	(12,312)(12,1212)
\path(792.000,1917.000)(762.000,2037.000)(732.000,1917.000)
\path(762,2037)(762,1212)
\path(762,2412)(762,1962)
\path(1662,1662)(1662,1212)
\path(1662,2412)(1662,1662)
\path(1632.000,1782.000)(1662.000,1662.000)(1692.000,1782.000)
\put(1012,1737){\makebox(0,0)[lb]{$\bar t$}}
\put(2087,1737){\makebox(0,0)[lb]{$ t$}}
\put(462,462){\makebox(0,0)[lb]{$\D(\nu^{1/p})$}}
\put(537,2562){\makebox(0,0)[lb]{$\Z(\gamma)$}}
\end{picture}
} }
         	\hspace{-1.9mm}
	\end{array}

$$
Now we calculate $\lift_p$ of the above sequence. Observe that both 
$\Z(\gamma)$ and $\nu^{1/p}$ are exponentials of series of connected 
diagrams with symmetric legs whose dashed graphs are not marked by any 
nontrivial beads. Thus,
one can check that $\lift_p$ is the closure of the following
diagram (there are $p$ copies displayed):
$$
\printname{tangle3}
	\setlength{\unitlength}{0.02\standardunitlength}
	\begin{array}{c}  \hspace{-1.7mm}
         	\raisebox{-8pt}{\begingroup\makeatletter\ifx\SetFigFont\undefined%
\gdef\SetFigFont#1#2#3#4#5{%
  \reset@font\fontsize{#1}{#2pt}%
  \fontfamily{#3}\fontseries{#4}\fontshape{#5}%
  \selectfont}%
\fi\endgroup%
{\renewcommand{\dashlinestretch}{30}
\begin{picture}(6699,6039)(0,-10)
\path(1662,312)(1662,12)
\path(762,312)(762,12)
\path(762,1212)(762,1512)
\path(1662,1212)(1662,1512)
\path(1662,2712)(1662,2412)
\path(762,2712)(762,2412)
\path(762,3612)(762,3312)
\path(1662,3612)(1662,3312)
\path(1662,4512)(1662,4812)
\path(762,4512)(762,4812)
\path(762,6012)(762,5712)
\path(1662,6012)(1662,5712)
\path(4962,3012)(4962,3312)
\path(5862,3012)(5862,3312)
\path(5862,4512)(5862,4212)
\path(4962,4512)(4962,4212)
\path(4962,2112)(4962,1812)
\path(5862,2112)(5862,1812)
\path(12,5712)(2412,5712)(2412,4812)
	(12,4812)(12,5712)
\path(12,4512)(2412,4512)(2412,3612)
	(12,3612)(12,4512)
\path(12,2412)(2412,2412)(2412,1512)
	(12,1512)(12,2412)
\path(12,1212)(2412,1212)(2412,312)
	(12,312)(12,1212)
\path(4212,3012)(6612,3012)(6612,2112)
	(4212,2112)(4212,3012)
\path(4287,4212)(6687,4212)(6687,3312)
	(4287,3312)(4287,4212)
\put(912,2862){\makebox(0,0)[lb]{$\dots$}}
\put(3237,2937){\makebox(0,0)[lb]{$=$}}
\put(387,537){\makebox(0,0)[lb]{$\D(\nu^{1/p})$}}
\put(387,1662){\makebox(0,0)[lb]{$\Z(\gamma)$}}
\put(387,3762){\makebox(0,0)[lb]{$\D(\nu^{1/p})$}}
\put(387,4962){\makebox(0,0)[lb]{$\Z(\gamma)$}}
\put(4587,3462){\makebox(0,0)[lb]{$\Z(\gamma)^p$}}
\put(4587,2262){\makebox(0,0)[lb]{$\D(\nu)$}}
\end{picture}
} }
         	\hspace{-1.9mm}
	\end{array}
=\Z(L^{(p)}).
$$
The proposition follows for $\Z$. The extension to the stated normalization
$\Zc$ is trivial.
\qed

The next proposition states that $\lift_p$ intertwines the integration
maps $\intrat$ and $\int$:

\begin{proposition}
\lbl{prop.liftint}
The following diagram commutes:
$$
%\divide\dgARROWLENGTH by2
\begin{diagram}
\node{\GA(\ostar_X,\Lloc)}
\arrow{e,t}{\intrat \,dX}\arrow{s,l}{\lift_p}
\node{\GA(\Lloc)}
\arrow{s,r}{\lift_p} \\
\node{\GA(\ostar_{X^{(p)}})}
\arrow{e,t}{\int \, dX^{(p)}}
\node{\GA(\phi)}
\end{diagram}
$$ 
\end{proposition}

Since $\intrat$, $\int$ and $\lift_p$ are partially defined maps (defined for 
$X$-integrable elements and for diagrams with nonsingular beads when evaluated
that complex $p$th roots of unity), the maps in the above diagram should be
restricted to the domain of definition of the maps, and the diagram then 
commutes, as the proof shows.

\begin{proof}
Consider a pair $(M,s)$ where $s$ is given by
$$
s=
\exp\left( \frac{1}{2} \sum_{i,j} 
%\eepic{wcoup}{0.02} 
\strutb{x_i}{x_j}{W_{ij}(t)}
\right) \sqcup R.
$$
If we write 
$$
 \lift_p(s) = \exp\left( \frac{1}{2} 
\sum_{i,j}\sum_{r=0}^{p-1}\sum_{s=0}^{p-1}
\st{x_i^{(r)}}{x_j^{(s)}}
W^{(p)}_{(i,r),(j,s)} 
%\eepic{strut2}{0.02}
\right)\sqcup R',
$$
then, observe that
\begin{eqnarray*}
\lift_p \left( 
%\eepic{wcoup}{0.02} 
\strutb{x_i}{x_j}{W_{ij}(t)}
\right) &=& \sum_{r=0}^{p-1}\sum_{s=0}^{p-1} 
%\eepic{strut2}{0.02} 
\st{x_i^{(r)}}{x_j^{(s)}} W^{(p)}_{(i,r),(j,s)}\\
%\,\,\,\,\,\,\, \text{ and } \,\,\,\,\,\,\,
%$$
%and
%$$
\lift_p(R) &=& R'.
\end{eqnarray*}

Recall the map $\ch_p: \Lloc^{(p)}\to\BC[t]/(t^p-1)$ of Equation 
\eqref{eq.chp}. 
It follows from the above that for any $r$ we have
$$
\ch_p(W_{ij}(t))=\sum_{s=0}^{p-1} W^{(p)}_{(i,r),(j,s)}t^{s-r}.
$$

We wish to determine $\ch_p(W_{ij}(t)^{-1})$, which we write as
$$
\ch_p(W_{ij}(t)^{-1})=\sum_{s=0}^{p-1} W^{(p)'}_{(i,r),(j,s)}t^{s-r}.
$$

Since $\delta_{ij}=\sum_k W_{ik} W^{-1}_{kj}$, we can solve for
$W^{(p)'}_{(i,r),(j,s)}$ in terms of $W^{(p)}_{(i,r),(j,s)}$ and obtain that

\begin{eqnarray*}
\lift_p\left( 
%\eepic{wcoupinv}{0.02} 
\strutb{x_i}{x_j}{W_{ij}^{-1}(t)}
\right) &=& 
\sum_{r=0}^{p-1}\sum_{s=0}^{p-1} 
\st{x_i^{(r)}}{x_j^{(s)}} (W^{(p)})^{-1}_{(i,r),(j,s)} .
%\eepic{strut2}{0.02}
\end{eqnarray*}

Observe further the following consequence of the ``state-sum'' definition of 
$\lift_p$: for diagrams $D_1,D_2$ in $\A(\star_X,\Lloc)$, we have that
$$
\lift_p \left( \la D_1, D_2 \ra_X \right)
=
\la \lift_p(D_1),\lift_p(D_2)\ra_{X^{(p)}}\ \in\ \A(\phi).
$$

Now, we can finish the proof of the proposition as follows:

\begin{eqnarray*}
\lift_p \left(\intrat dX \,\, (s) \right)
 & = & \lift_p \left( \left\la \exp
\left( -\frac{1}{2} \sum_{i,j} 
%\eepic{wcoupinv}{0.02}
\strutb{x_i}{x_j}{W_{ij}^{-1}(t)}
\right),
R \right\ra\right)_X \\
& = & \left\la \lift_p\left( \exp \left( -\frac{1}{2} \sum_{i,j} 
%\eepic{wcoupinv}{0.02}
\strutb{x_i}{x_j}{W_{ij}^{-1}(t)}
\right) \right) , \lift_p(R) \right\ra_{X^{(p)}}  \\
& = &
\left\la \exp \left( 
-\frac{1}{2} \sum_{i,j}\sum_{r=0}^{p-1}\sum_{s=0}^{p-1}  
\st{x_i^{(r)}}{x_j^{(s)}} (W^{(p)})^{-1}_{(i,r),(j,s)} 
%\eepic{strut2}{0.018}
\right) , R' \right\ra_{X^{(p)}} \\
& = & \int dX^{(p)} \,\,  \lift_p(s).
\end{eqnarray*}
\end{proof}

\begin{proof}(of Theorem \ref{thm.3})
It follows immediately from Propositions \ref{prop.fundlem} and 
\ref{prop.liftint}.
\end{proof}

\subsection{The connection of $\lift_p$ with $\bmod p$ residues.}
\lbl{sub.connection}

Rozansky in \cite{R} considered the following vector space $\A^R$
to lift the Kontsevich integral, 
$$
\A^R=\left(\oplus_\Ga   A_{\Ga,\loc}^R \cdot \Ga \right)/(\AS,\IHX)
$$
where the sum is over trivalent graphs $\Ga$ with oriented vertices and 
edges, and where 
$$
A_{\Ga,\loc}^R=(\BQ[\exp(H^1(\Ga,\BZ))]_{\loc})^{\Ga}
$$ 
is the $\Ga$-invariant subring of the 
(Cohn) localization of the group-ring $\BQ[\exp(H^1(\Ga,\BZ))]$ with 
respect to the ideal of elements that augment to $\pm 1$. We will think
of $A_{\Ga,\loc}^R$ as the coefficients by which a graph $\Ga$ is
multiplied.

Note that  $\BQ[\exp(H^1(\Ga,\BZ))]$ can be identified with the ring of 
Laurrent polynomials in $b_1(\Ga)$ variables, where $b_1(\Ga)$ is the first 
betti number of $\Ga$. Thus $\La_{\Ga,\loc}$ can be identified with the ring
of rational functions $p(s)/q(s)$ in $b_1(\Ga)$ variables $\{s\}$ for
polynomials $p$ and $q$ such that $q(1)=\pm 1$. Let $\La_{\Ga,\loc}^{(p)}$
denote the subring of $\La_{\Ga,\loc}$ that consists of functions $p(s)/q(s)$
as above such that $q$, evaluated at any complex $p$-th roots
of unity is nonzero.
In \cite{GR} (see also \cite{Kr2}) the authors considered a map:
$$
\Res_p: A_{\Ga,\loc}^{R,(p)} \to \BC
$$ defined by
$$
\Res_p \left( \frac{f(s)}{g(s)} \right)= p^{\chi(\Ga)}
\sum_{\w^p=1}\frac{f(\w)}{g(\w)}
$$
where the sum is over all $b_1(\Ga)$-tuples $(\w_1,\dots,w_{b_1(\Ga)})$
of complex $p$th root of unity and where $\chi(\Ga)$ is the Euler 
characteristic of $\Ga$.
This gives rise to a map $\Res_p: \A^R\to\A(\phi)$.

Similarly, we have that
$$
\A(\Lloc)= \left(\oplus_\Ga   A_{\Ga}(\Lloc) \cdot \Ga \right)/(
\mathrm{Relations})
$$
where $A_{\Ga}(\Lloc)$ is the $\Ga$-invariant subspace of the
vector space spanned by $\a : \mathrm{Edge}
(\Ga)\to \Lloc$ modulo the Relations of \cite[Figures 2,3]{GK2} which
include the $\AS,\IHX$ relations, multilinearity on the beads of the edges
and the Holonomy Relation.
%
%{\em Holonomy} Relation
%$$
%\psdraw{pushagain}{1.5in}
%$$
%where $r_i \in \Lloc$. 
An important difference between $\A^R$ and $\A(\Lloc)$ 
is the fact that $A_{\Ga,\loc}^R$ is an {\em algebra} 
whereas $A_{\Ga}(\Lloc)$ is only a {\em vector space}. Nevertheless, there 
is a map $\phi_{R,\Ga}: A_{\Ga}(\Lloc) \to A_{\Ga,\loc}^R$ defined by 
$$
\phi_{R,\Ga}(\a)=\frac{1}{\Aut(\Ga)} \sum_{\s \in \Aut(\Ga)} 
\prod_{e \in \mathrm{Edge(\Ga)}} \a_e(t_{\s(e)})
$$
where $\a=(\a_e(t)): \mathrm{Edge}(\Ga)\to\Lloc$. 
The maps $\phi_{R,\Ga}$ assemble together to define a map 
$\phi_R:\A(\Lloc)\to\A^R$.

For example, consider the trivalent graph $\Theta$ whose edges are
labeled by $\a,\b,\g \in \Lloc$ as shown below 
$$ 
\psdraw{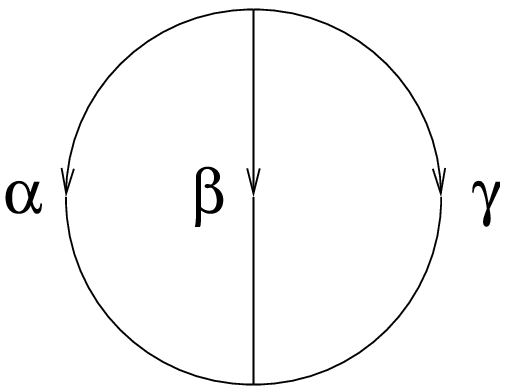}{1in} 
$$
with automorphism group
$\mathrm{Aut}(\Theta)=\text{Sym}_2 \times \text{Sym}_3$ that acts on 
the algebra of rational functions in three variables by permuting the
variables and by inverting all variables simultaneously. Then,
we have
$$
\phi_{R,\Theta}(\a,\b,\g)=\frac{1}{12} \sum_{\s \in 
\mathrm{Aut}(\Theta)} \a(t_{\s(1)})\b(t_{\s(2)})\g(t_{\s(3)}) 
\in \BQ(t_1,t_2,t_3).
$$

We finish by giving a promised relation between $\lift_p$ and $\Res_p$
for $p$-regular rational functions:

\begin{theorem}
\lbl{thm.liftres}
The following diagram commutes
$$
\divide\dgARROWLENGTH by2
\begin{diagram}
\node{\A(\Lloc^{(p)})}
\arrow[2]{e,t}{\phi_R}\arrow{se,r}{\lift_p}
\node[2]{\A^R} 
\arrow{sw,r}{\Res_p} \\
\node[2]{\A(\phi).}
\end{diagram}
$$
\end{theorem}

\begin{proof}
Using the properties of $\lift_p$ and $\Res_p$ it suffices to consider
only trivalent graphs $\Ga$ with edges decorated by elements in $\La$, 
and in fact only those graphs whose edges are decorated by powers of $t$.
Moreover, since both $\lift_p$ and $\Res_p$ satisfy the push relations,
it suffices to consider graphs whose edges along any forest are labeled by $1$.

Fix a trivalent graph $\Ga$ with ordered edges $e_i$ decorated by $\a=
(t^{m_1},\dots,t^{m_{3n}})$. We begin by giving a description of the algebra 
$A_{R,\loc}$
in terms of local coordinates as follows. Choose a maximal forest $T$ and 
assume, without loss of generality, that the edges of $\Ga\sminus T$ are 
$e_1,\dots,e_b$ where $b=b_1(\Ga)$. Each edge $e_i$ corresponds to a 
1-cocycle $x_i \in C^1(\Ga,\BQ)$. Since $H^1(\Ga,\BZ)=\Ker(C^1(\Ga,\BZ)\to 
C^0(\Ga,\BZ))$, it follows that $H^1(\Ga,\BZ)$ is a (free) abelian group
with generators $x_1,\dots,x_{3n}$ and  relations
$\sum_{j: \,\, v \in \partial e_j} \e_{j,v} x_j=0$ for all vertices $v$ of 
$\Ga$ and for appropriate local orientation signs $\e_{j,v} = \pm 1$.
It follows that 
$$
\BQ[H^1(\Ga,\BQ)]=\frac{\BQ[ t_{1}^{\pm 1}, \dots,t_{3n}^{\pm 1}]}{\left( 
 \prod_{j: \,\, v \in \partial e_j} t_j^{\e_{j,v}} =1
\text{ for } v \in \text{vertex}(\Ga) \right)}=
\BQ[ t_{1}^{\pm 1}, \dots,t_{b}^{\pm 1}],
$$
where $t_i=e^{x_i}$. This implies that $A_{\Ga,\loc}^R$ is a $\Ga$-invariant 
subalgebra of $\BQ(t_1, \dots,t_b)$.

Now, $\phi_{R,\Ga}(\a)$ is obtained by symmetrizing over $\Ga$-automorphisms
of the monomial $t_1^{m_1} \dots t_b^{m_b}$. We may assume that $m_i \in
\{0,\dots,p-1\}$ for all $i$. Thus,
$$
\Res_p(t_1^{m_1} \dots t_b^{m_b})=\frac{p^{b_0(\Ga)}}{p^b}
\left(\sum_{\w_1^p=1} \w_1^{m_1} \right)\dots
\left(\sum_{\w_b^p=1} \w_b^{m_b} \right)
=p^{b_0(\Ga)} 
\delta_{m_1,0} \dots \delta_{m_b,0}.
$$
 
On the other hand, an admissible $p$-coloring of $(\Ga,\a)$
necessarily assigns the same color to each connected component of $\Ga$
and then the consistency relations along the edges $e_i$ for $i=1,\dots,b$
show that an admissible coloring exists only if $m_i=0$, for $i=1,\dots,b$,
and in that case there is are $p$ admissible colorings for each connected
component of $\Ga$. Thus, the number of admissible $p$-colorings is
$p^b_0(\Ga)$.
 
After symmetrization over $\Ga$, the result follows.
\end{proof}

The reader is encouraged to compare the above proof with 
\cite[Lemmas 3.4.1, 3.4.2]{Kr2}.

\subsection{The degree $2$ part of $\Zrat$}
\lbl{sub.deg2}

In this section we prove Corollary \ref{cor.1}. 
The following lemma reformulates where $Q= \Zrat_2$ takes values. Consider the 
vector space
$$
\La_{\Th}= \otimes^3 \Lloc/\left( (f,g,h)=(tf,tg,th), \, \Aut(\Th)
\right)
$$
$\Aut(\Th)=\Sym_3 \times \Sym_2$ acts on $\otimes^3 \Lloc$ by permuting the 
three factors and by applying the involution of $\Lloc$ simultaneously
to all three factors.

\begin{lemma}
\lbl{lem.Q}
$Q$ takes values in $\La_{\Th} \cdot \Theta$
\end{lemma}

\begin{proof}
There are two trivalent graphs of degree $2$, namely $\Th$ and $\eyes$.
Label the three oriented edges of $\eyes$ $e_i$ for $i=1,2,3$ where $e_2$
is the label in the middle (nonloop) edge of $\eyes$. For $f,g,h \in \Lloc$,
let $\a_{\eyes}(f,g,h) \in \a_{\eyes}(\Lloc) \cdot \eyes$ denote the corresponding
element.

For $p,q \in \La$, $f,h \in \Lloc$, we write $q=\sum_k a_k t^k$ and compute

\begin{xalignat*}{2} 
\a_{\eyes}(f,p,h) &= \a(f, (p/q).q, h) & \\
&= \sum_k \a(f,(p/q) a_k t^k, h) & \text{by Multilinearity} \\
&= \sum_k \a_(f, p/q, a_k t^k h t^{-k}) & \text{by the Holonomy Relation} \\
&= \a(f,p/q, q(1) h) &
\end{xalignat*}

Thus, $\a_{\eyes}(\Lloc)$ is spanned by $\a(f,p,h)$ for $f,p,h$ as above.
Applying the above reasoning once again, it follows that $\a_{\eyes}(\Lloc)$
is spanned by $\a(f,1,h)$ for $f,h$ as above.

Applying the $\IHX$ relation
$$
\psdraw{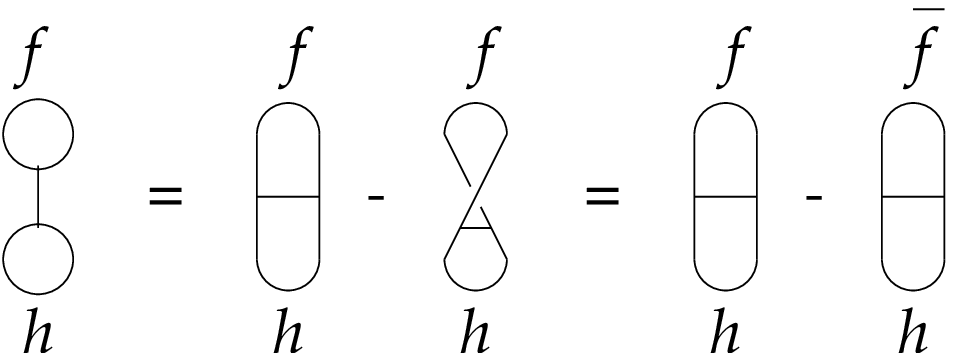}{2.5in}
$$
it follows that the natural map $\La_{\Th} \to \A_2(\Lloc)$ is onto.
It is easy to see that it is also 1-1, thus a vector space isomorphism.
\end{proof}

\begin{remark}
\lbl{rem.deg2}
In fact, one can show that $Q$ takes values in the abelian subgroup
$\La_{\Th,\BZ}$ of $\La_{\Th}$ generated by $\otimes^3 \Lz$.
\end{remark}

\begin{proof}(of Corollary \ref{cor.1})
Consider the degree $2$ part in the 
Equation of Theorem \ref{thm.1}. On the one hand, we have
$Z_2=1/2\lb \cdot \Th$ (see \cite[Section 5.2]{LMO}) and 
on the other hand, it follows by definition and Lemma \ref{lem.Q} that
$\Zrat_2=1/6 Q \cdot \Th$. Theorem \ref{thm.liftres} which compares liftings
ad residues concludes first part of the corollary.

For the second part, observe that $Q$ is a rational function on 
$S^1 \times S^1$, which is regular when evaluated at complex roots of
unity. Furthermore, by definition of $\Res_p$, it follows that
$$
\frac{1}{p}\Res_p^{t_1,t_2,t_3}Q(M,K)=\frac{1}{p^2} 
\sum_{\w_1^p=\w_2^p=1}Q(M,K)(\w_1,\w_2,(\w_1 \w_2)^{-1})
$$
is the average of $Q(M,K)$ on $S^1 \times S^1$ (evaluated at pairs of
complex $p$th
roots of unity) and converges to $\int_{S^1 \times S^1} Q(M,K)(s) d\mu(s)$.
This concludes the proof of Corollary \ref{cor.1}.
\end{proof}

%We finish with the following question, which seems obvious, but is yet 
%unknown to us:

%\begin{question}
%\lbl{que.1}
%Is the map $\hair: \A(\Lloc)\to\A(\star)$ one-to-one?
%\end{question}

\section{Remembering the knot}
\lbl{sec.remember}

In this section we will briefly discuss an extension of Theorem \ref{thm.1}
for invariants of cyclic branched covers in the presence of the lift
of the branch locus. 

We begin by noting that the rational invariant $\Zrat$ can be extended to
an invariant of pairs $(M,K)$ of {\em null homologous knots $K$ in \qhs s
$M$}, \cite{GK2}. The extended invariant (which we will denote by the same
name), takes values in $\GAz(\Lloc)=\B(\Lz\to\BQ) \times \GA(\Lloc)$. 
In this section, we will work in this generality.

Consider a pair $(M,K)$ of a null homologous knot $K$ in a \qhs\ $M$, 
and the corresponding cyclic branched covers $\S^p_{(M,K)}$. The preimage of 
$K$ in $\S^p_{(M,K)}$ is a knot $\kbr$, which we claim is null homologous. 
Indeed, we can construct the branched 
coverings by cutting $M-K$ along a Seifert surface of $K$ and gluing several
copies side by side. This implies that a Seifert surface of $K$ in $M$
lifts to a Seifert surface of $\kbr$ in $\S^p_{(M,K)}$. 

If we wish, we may think of $\kbr$ as a $0$-framed knot in $\S^p_{(M,K)}$
(where a $0$-framing is obtained by a parallel of $\kbr$ along a Seifert
surface, and is independent of the Seifert surface chosen).

We now consider the rational invariant $\Zrat(\S^p_{(M,K)},\kbr)$
of a {\em $p$-regular pair} $(M,K)$, that is a pair such that $M$ and 
$\S^p_{(M,K)}$ are \qhs s and $K$ is null homologous in $M$. 
For the rational version of the lift map
$$
\liftrat_p: \GAz(\Lloc) \to \GAz(\Lloc).
$$
defined below, we have the following improved version of Theorem \ref{thm.1}:

\begin{theorem}
\lbl{thm.11}
For all $p$ and $p$-regular pairs $(M,K)$, we have
$$\Zrat(\Sp_{(M,K)},\kbr)=e^{\s_p(M,K)\Theta/16} 
\liftrat_p \circ \twistrat_{\a_p} \circ \Zrat(M,K) \in  
\GAz(\Lloc).
$$
where $\a_p=\nu^{-(p-1)/p}$, $\nu=Z(S^3,\mathrm{unknot})$.
\end{theorem}

The meaning of multiplying elements $(M,s) \in \GAz(\Lloc)$ by elements $a \in
\A(\phi)$ is as follows: $a\cdot (M,s)=(M, a \sqcup s)$.

\begin{remark}
Evaluating $\GAz(\Lloc)\to\GA(\phi)$ at $t=1$ corresponds to
forgetting the knot $\kbr$, thus the above theorem is an improved version
of Theorem \ref{thm.1}.

The proof of Theorem \ref{thm.11}, which is left as an exercise,
follows the same lines as the proof
of Theorem \ref{thm.1} using properties of the $\liftrat_p$ map rather
than properties of the $\lift_p$ map.
\end{remark}

In the remaining section, we introduce the map $\liftrat_p$ which is an 
enhancement of the map $\lift_p$ of Section \ref{sec.lifting}.
We start by defining a map
$$
\liftrat_p: \A(\tline_X,\La)\to \A(\tline_{X^{(p)}},\La).
$$
This map is defined in exactly the same way as the map $\lift_p$ of Lemma
\ref{lem.liftp1}, except that instead of forgetting all labels as the last
step, we do the following replacement:
$$
\printname{treplace}
	\setlength{\unitlength}{0.03\standardunitlength}
	\begin{array}{c}  \hspace{-1.7mm}
         	\raisebox{-8pt}{\begingroup\makeatletter\ifx\SetFigFont\undefined%
\gdef\SetFigFont#1#2#3#4#5{%
  \reset@font\fontsize{#1}{#2pt}%
  \fontfamily{#3}\fontseries{#4}\fontshape{#5}%
  \selectfont}%
\fi\endgroup%
{\renewcommand{\dashlinestretch}{30}
\begin{picture}(4708,1539)(0,-10)
\put(4337,762){$\bullet$}%\ellipse{150}{150}}
\put(-27,762){$\bullet$}%\ellipse{150}{150}}
\dashline{60.000}(87,1512)(87,12)
\path(12,1437)(87,1512)(162,1437)
\dashline{60.000}(2337,1512)(2337,12)
\path(2262,1437)(2337,1512)(2412,1437)
\path(837,762)(1737,762)
\path(1617.000,732.000)(1737.000,762.000)(1617.000,792.000)
\dashline{60.000}(4437,1512)(4437,12)
\path(4362,1437)(4437,1512)(4512,1437)
\put(312,87){\makebox(0,0)[lb]{$a$}}
\put(312,1287){\makebox(0,0)[lb]{$a+1$}}
\put(312,687){\makebox(0,0)[lb]{$t$}}
\put(4662,687){\makebox(0,0)[lb]{$t$}}
\put(3087,687){\makebox(0,0)[lb]{$\text{or}$}}
\end{picture}
} }
         	\hspace{-1.9mm}
	\end{array}
 
$$
depending on $a \neq p-1$ or $a=p-1$.
As in Section \ref{sub.deflift}, this leads to a well-defined map
$$
\liftrat_p: \GA(\ostar_X,\La)\to \GA(\ostar_{X^{(p)}},\La).
$$
The next step is to extend this to a map of diagrams with rational beads
in $\Lloc^{(p)}$. The following lemma considers elements of the ring
$\Lloc^{(p)}$.

\begin{lemma}
\lbl{lem.llocp}
Every $r(t) \in \Lloc^{(p)}$ can be written in the form
$
r(t)=p(t)/q(t^p)
$
where $p(t),q(t) \in \La \otimes \BC$.
\end{lemma}

\begin{proof}
Using a partial fraction expansion of the denominator of $r(t)$, it
suffices to assume that $r(t)=1/(t-a)^k$ for some $k \geq 1$. In that case, 
we have
$$
\frac{1}{t-a} = \frac{ \prod_{i=1}^{p-1} (t- a \omega^i)}{t^p-a^p}
$$
where $\omega=\exp(2 \pi i/p)$.
\end{proof}
Now, we can introduce the definition of $\liftrat_p$ for diagrams with
labels in $\Lloc^{(p)}$. Consider such a diagram $D$, and replace each
bead $r(t)$ by a product of beads $p(t) \, 1/q(t^p)$ using Lemma 
\ref{lem.llocp}.
Now, consider the diagrams obtained by $p$-admissible colorings of the lift 
$\pi_p^{-1}(D)$, that is colorings of the lift that satisfy the following 
conditions:
$$
\printname{treplace4}
	\setlength{\unitlength}{0.03\standardunitlength}
	\begin{array}{c}  \hspace{-1.7mm}
         	\raisebox{-8pt}{\begingroup\makeatletter\ifx\SetFigFont\undefined%
\gdef\SetFigFont#1#2#3#4#5{%
  \reset@font\fontsize{#1}{#2pt}%
  \fontfamily{#3}\fontseries{#4}\fontshape{#5}%
  \selectfont}%
\fi\endgroup%
{\renewcommand{\dashlinestretch}{30}
\begin{picture}(8321,1539)(0,-10)
\put(5462,762){$\bullet$}%\ellipse{150}{150}}
\dashline{60.000}(762,1512)(762,12)
\dashline{60.000}(2412,1512)(3012,762)(3612,1512)
\dashline{60.000}(3012,12)(3012,762)
\dashline{60.000}(5562,1512)(5562,12)
\path(12,12)(1512,12)
\path(1392.000,-18.000)(1512.000,12.000)(1392.000,42.000)
\path(5487,1437)(5562,1512)(5637,1437)
\put(462,1062){\makebox(0,0)[lb]{$a$}}
\put(2337,1062){\makebox(0,0)[lb]{$a$}}
\put(3237,87){\makebox(0,0)[lb]{$a$}}
\put(5787,87){\makebox(0,0)[lb]{$a$}}
\put(5787,1287){\makebox(0,0)[lb]{$a+1$}}
\put(1212,237){\makebox(0,0)[lb]{$a$}}
\put(3537,1062){\makebox(0,0)[lb]{$a$}}
\put(237,237){\makebox(0,0)[lb]{$a$}}
\put(5787,687){\makebox(0,0)[lb]{$t$}}
\put(7412,762){$\bullet$}%\ellipse{150}{150}}
\dashline{60.000}(7512,1512)(7512,12)
\path(7437,1437)(7512,1512)(7587,1437)
\put(7737,687){\makebox(0,0)[lb]{$1/q(t^p)$}}
\put(7812,1287){\makebox(0,0)[lb]{$a$}}
\put(7812,12){\makebox(0,0)[lb]{$a$}}
\end{picture}
} }
         	\hspace{-1.9mm}
	\end{array}

$$
Finally forget the beads of the edges, as follows:
$$
\printname{treplace6}
	\setlength{\unitlength}{0.03\standardunitlength}
	\begin{array}{c}  \hspace{-1.7mm}
         	\raisebox{-8pt}{\begingroup\makeatletter\ifx\SetFigFont\undefined%
\gdef\SetFigFont#1#2#3#4#5{%
  \reset@font\fontsize{#1}{#2pt}%
  \fontfamily{#3}\fontseries{#4}\fontshape{#5}%
  \selectfont}%
\fi\endgroup%
{\renewcommand{\dashlinestretch}{30}
\begin{picture}(11732,1581)(0,-10)
\put(83,804){$\bullet$}%\ellipse{150}{150}}
\put(4433,804){$\bullet$}%\ellipse{150}{150}}
\put(6533,804){$\bullet$}%\ellipse{150}{150}}
\put(8633,804){$\bullet$}%\ellipse{150}{150}}
\put(11033,804){$\bullet$}%\ellipse{150}{150}}
\dashline{60.000}(158,1554)(158,54)
\path(83,1479)(158,1554)(233,1479)
\dashline{60.000}(2408,1554)(2408,54)
\path(2333,1479)(2408,1554)(2483,1479)
\path(908,804)(1808,804)
\path(1688.000,774.000)(1808.000,804.000)(1688.000,834.000)
\dashline{60.000}(4508,1554)(4508,54)
\path(4433,1479)(4508,1554)(4583,1479)
\path(5108,804)(6008,804)
\path(5888.000,774.000)(6008.000,804.000)(5888.000,834.000)
\dashline{60.000}(6608,1554)(6608,54)
\path(6533,1479)(6608,1554)(6683,1479)
\dashline{60.000}(8708,1554)(8708,54)
\path(8633,1479)(8708,1554)(8783,1479)
\path(9683,804)(10583,804)
\path(10463.000,774.000)(10583.000,804.000)(10463.000,834.000)
\dashline{60.000}(11108,1554)(11108,54)
\path(11033,1479)(11108,1554)(11183,1479)
\put(383,129){\makebox(0,0)[lb]{$a \neq p-1$}}
\put(383,1329){\makebox(0,0)[lb]{$a+1$}}
\put(383,729){\makebox(0,0)[lb]{$t$}}
\put(4733,729){\makebox(0,0)[lb]{$t$}}
\put(6833,729){\makebox(0,0)[lb]{$t$}}
\put(4733,54){\makebox(0,0)[lb]{$p-1$}}
\put(4733,1329){\makebox(0,0)[lb]{$0$}}
\put(8933,1329){\makebox(0,0)[lb]{$b$}}
\put(8933,54){\makebox(0,0)[lb]{$b$}}
\put(8933,729){\makebox(0,0)[lb]{$1/q(t^p)$}}
\put(11333,729){\makebox(0,0)[lb]{$1/q(t)$}}
\end{picture}
} }
         	\hspace{-1.9mm}
	\end{array}
 \hspace{0.5cm} .
$$
$\liftrat_p(D)$ is defined to be the resulting combination of diagrams.
We leave as an exercise to show that this is well-defined, independent
of the quotient used in Lemma \ref{lem.llocp} above.

\begin{remark}
\lbl{rem.liftrat}
The map $\liftrat_p: \GA(\ostar_X,\Lloc^{(p)})\to \GA(\ostar_X,\Lloc)$
is an algebra map, using thr disjoint union multiplication.
\end{remark}

Finally, we define
$$
\liftrat_p: \GAz(\ostar_X,\Lloc^{(p)})\to \GAz(\ostar_{X^{(p)}},\Lloc)
$$
by 
$$
\liftrat_p(M(t), s)=(M(t \to  T_t^{(p)}), \liftrat_p(s)),
$$
where $T_t^{(p)}$ is the $p$ by $p$ matrix (given by example for $p=4$)
\begin{equation}
\lbl{eq.Tpp}
T_t^{(p)}=\left[
\begin{array}{llll}
0 & 1 & 0 & 0 \\
0 & 0 & 1 & 0 \\
0 & 0 & 0 & 1 \\
t & 0 & 0 & 0 
\end{array}
\right] .
\end{equation}

The substitution of Equation \eqref{eq.Tpp} is motivated from the 
combinatorics of lifting struts (the analogue of Proposition 
\ref{prop.liftint} for the $\liftrat_p$ map), but also from the following
lemma from algebraic topology, that was communicated to us by J. Levine,
and improved our understanding:

\begin{lemma}
\lbl{lem.Tpp}
Consider a null homotopic link $L$ in a standard solid torus $ST$, with
equivariant linking matrix $A(t)$ and its lift $L^{(p)}$ in $ST$ under the 
$p$-fold covering map $\pi_p: ST\to ST$. Then, $L^{(p)}$ is null homotopic in
$ST$ with equivariant linking matrix given by $A(t \to T_t^{(p)})$.
\end{lemma}

\begin{proof}
Consider the commutative diagram
$$
%\divide\dgARROWLENGTH by2
\begin{diagram}
\node{\widetilde{ST}}
\arrow{e,t}{\pi'}\arrow{se,t}{\pi} 
\node{ST} 
\arrow{s,r}{\pi_p} \\
\node[2]{ST}
\end{diagram}
$$ 
where $\pi$ and $\pi'$ are universal covering maps. Since $\pi_p$
is 1-1 on fundamental groups, it follows that $L^{(p)}$ is null homotopic
in $ST$. Choose representatives $L'_i$ of the components $L_i$
of $L$ in the unversal cover $\widetilde{ST}$, for $i=1,\dots,l$
where $l$ is the number of components of $L$. Then, 
$$
A_{ij}(t)=\sum_{k=0}^\infty \lk( L'_i, t^k L'_j)t^k.
$$
On the other hand, 
$\{ t^r L'_i \}$ is a choice of representatives of the lifts of $L^{(p)}$
to $\widetilde{ST}$, for $r=0,\dots,p-1$, and $l=1,\dots,l$.
Furthermore, if $(B_{ij,rs}(t))$ is the equivariant linking matrix of 
$L^{(p)}$, we have
$$
B(t)_{ij,rs}(t)=\sum_{k=0}^\infty \lk( t^r  L'_i, t^{k+j} L'_j)t^k.
$$
It follows that if we collect all powers of $t$ modulo $p$ in Laurrent
polymomials $a_{ijk}$ such that
$$
t^{r-s} A_{ij}(t)=\sum_{k=0}^{p-1} a_{ij,rs,k}(t^p)
$$
(for $r,s =0,\dots,p-1$), then
$$
B_{ij,rs}(t)= a_{ij,rs,0}.
$$
Writing this in matrix form, gives the result.
\end{proof}

We end this section with a comment regarding the commutativity of
$\twistrat$ and $\liftrat_p$ as endomorphisms of $\GA(\Lloc^{(p)})$:

\begin{lemma}
\lbl{lem.commute}
For $\a \in \GA(\star)$, we have
$$
\liftrat_p \circ \twistrat_{\a} = \twistrat_{\a^p} \circ \liftrat_p .
$$
\end{lemma}

\section{The wheeled invariants}
\lbl{sec.wheel}

The goal of this independent section is to discuss the relation between
twisting and wheeling of diagrams and, as an application, to give an 
alternative version of Theorem \ref{thm.1} in terms of the wheeled rational
invariant $\Zratw$ introduced below.

Recall the Wheeling and Unwheeling maps from Section \ref{sub.magic}.

\begin{lemma}
\lbl{lem.twistwheel}
For $x \in \A(\star)$, we have
\begin{eqnarray*}
\twist_{\Om}(x) &=& \la \Om,\Om \ra^{-1} \,  x^{\wheel^{-1}} \\
\twist_{\Om^{-1}_\#}(x) &=& \la \Om,\Om \ra \, x^{\wheel}
\end{eqnarray*}
where the notation $\Om^{-1}_\#$ means the inverse of $\Om \in \A(\star)$
using the $\#$ multiplication (rather than the disjoint union multiplication).
\end{lemma}

Note that $\chi(\Om^{r}_\#)=\nu^{r}$, for all $r \in \BQ$, by notation.

\begin{proof}
The first identity follows from Lemma \ref{lem.twist}(a) using the identity
$$
\homega^{-1}(\Om)=\la \Om,\Om \ra^{-1} \Om
$$
of \cite[Proposition 3.3 and Corollary 3.5]{BL}.

The second identity follows from the first, after inverting the operators
involved. Specifically, Lemma \ref{lem.twist}(b) implies that
\begin{eqnarray*}
y & = & \twist_{1}(y) \\
& = & \twist_{\Om^{-1}_\#}(\twist_\Om(y)) \\
& = & \twist_{\Om^{-1}_\#}(\la \Om,\Om \ra^{-1} \, y^{\wheel^{-1}}) \\
& = & \la \Om,\Om \ra^{-1} \, \twist_{\Om^{-1}_\#}(y^{\wheel^{-1}})
\end{eqnarray*}
Setting $x =y^{\wheel^{-1}}$, we have that $y=x^{\wheel}$ and the above
implies that
$$
x^{\wheel}=\la \Om,\Om \ra^{-1} \, \twist_{\Om^{-1}_\#}(x).
$$
\end{proof}

The {\em wheeled invariant} $\Zw$ is defined by wheeling the $\Z$ invariant 
of each of the component of a link. 
Although $\Zw$ is an invariant of links equivalent to the $Z$ invariant, 
in many cases the $\Zw$ invariant behaves in a more natural way, as was 
explained in \cite{BL}. Similarly, we define the 
{\em wheeled rational invariant} $\Zratw$ by 
$$
\Zratw(M,K)=\twistrat_{\Om^{-1}_{\#}} \circ \Zrat(M,K) \in \GAz(\Lloc).
$$
The naming of $\Zratw$ is justified by the following equation
$$
\hairnu \circ \Zratw(M,K) = \la \Om,\Om \ra \, \Zw(M,K) \in \A(\star)
$$
which follows from Corollary \ref{cor.ZtoZ} (with $\a=\Om^{-1}_{\#}$)
and Lemma \ref{lem.twistwheel}.

The rational wheeled invariant $\Zratw$ behaves in some ways more naturally
than the $\Zrat$ invariant. A support of this belief is the 
following version of Theorem \ref{thm.11}:

\begin{theorem}
\lbl{thm.1alt}
For all $p$ and $p$-regular pairs $(M,K)$ we have
$$\Zratw(\Sp_{(M,K)},\kbr)=e^{\s_p(M,K)\Theta/16} \,
\lift_p \circ \Zratw(M,K) \in \GAz(\Lloc).
$$
\end{theorem}

The proof uses the same formal calculation that proves Theorem \ref{thm.1},
together with the following version of Theorem \ref{thm.3}:

\begin{theorem}
\lbl{thm.3alt}
With the notation of Theorem \ref{thm.3}, we have
$$
\lift_p \left(\intrat dX \,\, \Zratcw(L) \right) 
= \int dX^{(p)} \,\,  \Zcw(L^{(p)}).
$$
\end{theorem}

The proof of Theorem \ref{thm.3alt} follows from the proof of 
Proposition \ref{prop.fundlem}. %and \ref{prop.twistwheel} above.

\appendix

\section{Diagrammatic calculus}
\lbl{sec.repeat}

In this section we finish the proof of Theorem \ref{prop.ztwist}
using the identities and using the notation of \cite[Appendices A-E]{GK2}.
The rest of the proof uses the function-theory
properties of the $\intrat$-integration \cite[Appendix A-E]{GK2}.
These properties are expressed in terms of
the combinatorics of gluings of legs of diagrams.
$\intrat$-integration is a diagrammatic Formal Gaussian Integration that
mimics closely the Feynmann diagram expansion of perturbative quantum
field theory. Keep this in mind particularly with manipulations below
called the ``$\delta$-function trick'', ``integration by parts lemma'' and 
``completing the square''. The uninitiated reader may consult 
\cite[Part I,II]{A} for examples and motivation of the combinatorial calculus
and also \cite{GK2,BLT,BL}. We will follow the notation of \cite{A,GK2} here.

We focus on the term $\intrat dX \left(\varphi (\Zratc(L)) \right).$
Let us assume that the canonical decomposition of $\Zratc(L)$ is
$$
\Zratc(L) = \exp\left( \frac{1}{2} \sum 
%\gau{W_{ij}(t)}{x_j}{x_i}{\ \ \ \ \ }
\strutb{x_j}{x_i}{W_{ij}}   
\right)
\sqcup R,
$$
suppressing summation indices. 
We perform a standard move (the ``$\delta$-function trick'') to write this as:
$$
\left\la R(y), \exp\left( \frac{1}{2} \sum 
%\gau{W_{ij}(t)}{x_j}{x_i}{\ \ \ \ \ }
\strutb{x_j}{x_i}{W_{ij}}
+ 
%\gau{}{y_i}{x_i}{} 
\strutb{y_j}{x_i}{}   
\right) \right\ra_Y,
$$
where $Y$ is a set in 1-1 correspondence with $X$. Continuing,
\begin{eqnarray*}
\intrat dX \left(\varphi (\Zratc(L)) \right)
& = &
\left\la
\varphi\left(
R(y)\right),
\intrat dX \left( \exp\left( \frac{1}{2} \sum
\varphi
\left(
%\gau{W_{ij}(t)}{x_j}{x_i}{\ \ \ \ \ }
\strutb{x_j}{x_i}{W_{ij}}
\right)
+ 
%\gau{}{y_i}{x_i}{}
\strutb{y_i}{x_i}{}
\right) \right)
\right\ra_Y .
\end{eqnarray*}

\noindent
The ``integration by parts lemma'' \cite{GK2} implies that
\begin{multline*}
\intrat dX \left( 
\exp\left( \frac{1}{2} \sum
\varphi
\left(
\strutb{x_j}{x_i}{W_{ij}}
\right)
+ 
\strutb{y_j}{x_i}{}
\right) \right)  = \\ 
\intrat dX \left( 
\left(
\exp\left(
\sum
\varphi\left(
- 
\strutb{x_j}{y_i}{W^{-1}_{ij}}
\right)
\right) \right)\flat_{X}
\exp\left( \frac{1}{2} \sum
\varphi
\left(
\strutb{x_j}{x_i}{W_{ij}}   
\right)
+ 
\strutb{y_i}{x_i}{}
\right) \right) .
\end{multline*}

\noindent
``Completing the square'' implies that the above equals to:

$$
\exp\left( -\frac{1}{2} \sum \varphi \left(
%\gau{W_{ij}^{-1}(t)}{y_j}{y_i}
\strutb{y_j}{y_i}{W_{ij}^{-1}}
\right) \right)
\,\,\,\,\,
\intrat dX \left(
\exp\left(
\frac{1}{2}
\sum
\varphi \left(
%\gau{W_{ij}(t)}{x_j}{x_i}{\ \ \ \ \ }
\strutb{x_j}{x_i}{W_{ij}}
\right)
\right)
\right).
$$

\noindent
Returning to the expression in question:

\begin{eqnarray*}
\intrat dX \left(\varphi (\Zratc(L)) \right)
& = & \intrat dX \left(
\exp\left(
\frac{1}{2}
\sum
\varphi \left(
\strutb{x_j}{x_i}{W_{ij}}
\right)
\right)
\right) \\
& &  \sqcup 
\left\la
\varphi \left(
R(y)\right),
\exp\left( -\frac{1}{2} \sum \varphi \left(
\strutb{y_j}{y_i}{W_{ij}^{-1}}
\right) \right) \right\ra_Y. 
\end{eqnarray*}

\noindent
The second factor equals to $\varphi \left( \Zrat(M,K) \right)$.
The first factor contains only sums of disjoint union of wheels.
We can repeat the arguments which lead to the 
proof of the of the Wheels Identity in this case, \cite[Appendix E]{GK2}.

\begin{eqnarray*}
\lefteqn{
\intrat dX \left(
\exp\left(
\frac{1}{2}
\sum
\varphi \left(
%\gau{W_{ij}(t)}{x_j}{x_i}{\ \ \ \ \ }
\strutb{x_j}{x_i}{W_{ij}}
\right)
\right)
\right)} \\
& = &
\intrat dX \left( 
\exp \left(
\frac{1}{2}
\sum
\left( 
%\gau{W_{ij}(t)}{x_j}{x_i}{\ \ \ \ \ }
\strutb{x_j}{x_i}{W_{ij}}
+
\left(
\varphi \left(
%\gau{W_{ij}(t)}{x_j}{x_i}{\ \ \ \ \ }
\strutb{x_j}{x_i}{W_{ij}}
\right)
-
%\gau{W_{ij}(t)}{x_j}{x_i}{\ \ \ \ \ }
\strutb{x_j}{x_i}{W_{ij}}
\right) \right) \right) \right) \\
& = &
\left\la
\exp \left(
-\frac{1}{2}
\sum
%\gau{W_{ij}^{-1}(t)}{ x_j}{ x_i}{\ \ \ \ \ }
\strutb{ x_j}{ x_i}{W_{ij}^{-1}}
\right),
\exp \left( \frac{1}{2} \left( 
\varphi \left(
%\gau{W_{ij}(t)}{x_j}{x_i}{\ \ \ \ \ }
\strutb{x_j}{x_i}{W_{ij}}
\right)
-
%\gau{W_{ij}(t)}{x_j}{x_i}{\ \ \ \ \ }
\strutb{x_j}{x_i}{W_{ij}}
\right)\right)
\right\ra_X \\
%& & \mbox{{\bf figure out}}\\
& = &
\exp \left( -\frac{1}{2} \tr\log\left( W^{-1} \varphi(W) \right) 
\right)  .
% \\
% & = &
% \varphi \exp \left( -\frac{1}{2} \tr\log
% \left( \varphi(W)^{-1} W \right) 
% \right) . 
\end{eqnarray*}

\qed

% {\bf Todo list:

% \begin{itemize}
%\item
%ambiguous notation: $\Zrat(L), \Zrat(M,K)$.
%\item
%Do we need group-like in the definition of $\twist_\a$ and $\twistrat_\a$?
%\item
%be careful with factors of $\Om$ in twisting.
%\item
%Maybe twisting simplifies the notation once you put the $\Om$'s in the 
%right place?
%\item
%Missing the proof that $\lift_p$ respects the basing relations.
%\item
%Extend $\twistrat$ on $\GAz$. Use it to clean out the statement of Theorem 6.
%\item
%Make sure that $\phi_{t\to te^h}:\GAz(\star_X,\Lloc)\to\GAz(\star_{X \cup 
%h},\Lloc)$ 
%respects the push and the group-like basing relations.
%\item
%Group-like elements in the definition of $\lift_p$?
%\item
%Figure out the growth rates of the Casson invariant as functions of $p$.
%Is it exponential growth?
%\item
%Make a comment that $\twistrat$ and $\lift_p$ are intertwiners of integration
%theories. Intertwiners are useful objects.
%\item
%Replace $\nu$ by $\Om$.
%\item
%In section \ref{sub.magic} mention that $\hat{\a}$ preserves link relations
%and can be defined in the presence of other colors.
%\item
%Decide to write $A_x$ versus $A(x)$.
%\item
%Erase the notation $\Z[\a]$ and $\Zrat[\a]$ and replace it with $\twist_\a$
%and $\twistrat_\a$.
% \item
% What is the importance of $\twistrat$?
% \end{itemize}
% }

\ifx\undefined\bysame
	\newcommand{\bysame}{\leavevmode\hbox
to3em{\hrulefill}\,}
\fi

\end{document}